\documentclass[12pt,english]{article}
\usepackage{fullpage}
\usepackage{babel}
\usepackage{exscale,xspace}
\usepackage{amsfonts,amssymb}

\usepackage[
  pdfauthor={Bernd C. Kellner},
  pdfkeywords={Bernoulli numbers, irregular pairs of higher order, index of irregularity,
    Kummer congruences, p-adic zeta function, Riemann zeta function,
    Iwasawa theory, Teichmuller character},
  pdftitle={On irregular prime power divisors of the Bernoulli numbers},
  colorlinks,bookmarksnumbered]{hyperref}

\newcommand{\binom}[2]{ {#1 \choose #2} }
\newcommand{\pdiv}{\mid}
\newcommand{\pdiveq}{\mid\mid}
\newcommand{\notdiv}{\nmid}
\newcommand{\mods}[1]{({\rm mod\ } #1)}
\newcommand{\modp}[1]{\ ({\rm mod\ } #1)}
\newcommand{\qed}{\parbox{0cm}{}\hspace*{\fill} $\Box$}

\def\ord{\mathop{\rm ord}\nolimits}
\def\sign{\mathop{\rm sign}\nolimits}
\def\denom{\mathop{\rm denom}\nolimits}
\def\numer{\mathop{\rm numer}\nolimits}

\def\modres{\,\mathop{\rm mod}\,}

\newcommand{\QQ}{\mathbb{Q}}
\newcommand{\ZZ}{\mathbb{Z}}
\newcommand{\NN}{\mathbb{N}}

\newcommand{\PP}{\mathbb{P}}

\newcommand{\refeqn}[1]{(\ref{#1})}

\newcommand{\BN}{\widehat{B}}
\newcommand{\IRR}{\Psi^{\rm irr}}
\newcommand{\IRP}{\widehat{\Psi}^{\rm irr}}
\newcommand{\eulerphi}{\varphi}
\newcommand{\Bnum}{\mathit{\Lambda}}
\newcommand{\Bden}{V}
\newcommand{\Dop}{\mathcal{D}}

\renewcommand{\theequation}{\arabic{section}.\arabic{equation}}
\newcommand{\settheequation}[1]{\renewcommand{\theequation}{#1}}
\newcommand{\restoretheequation}{\renewcommand{\theequation}%
{\arabic{section}.\arabic{equation}}\addtocounter{equation}{-1}}

\newtheorem{prop}{Proposition}[section]
\newtheorem{theorem}[prop]{Theorem}
\newtheorem{corl}[prop]{Corollary}
\newtheorem{lemma}[prop]{Lemma}
\newtheorem{conj}[prop]{Conjecture}

\newtheorem{thdefin}[prop]{Definition}
\newtheorem{thremark}[prop]{Remark}
\newtheorem{thtable}[prop]{Table}
\newtheorem{thdiag}[prop]{Diagram}
\newtheorem{thexamp}[prop]{Example}

\newenvironment{defin}{\begin{thdefin}\rm}{\end{thdefin}}
\newenvironment{remark}{\begin{thremark}\rm}{\end{thremark}}
\newenvironment{tablenv}{\begin{thtable}\rm}{\end{thtable}}
\newenvironment{diagenv}{\begin{thdiag}\rm}{\end{thdiag}}
\newenvironment{expenv}{\begin{thexamp}\rm}{\end{thexamp}}

\newenvironment{proof}{\medskip\noindent%
\textsc{Proof.}}{\qed\medskip}
\newenvironment{proofof}[1]{\medskip\noindent%
\textsc{Proof of #1.}}{\qed\medskip}

\newenvironment{equationname}[1]{\settheequation{#1}\begin{equation}}%
{\end{equation}\restoretheequation}

\begin{document}

\title{On irregular prime power divisors of the\\ Bernoulli numbers}
\author{Bernd C. Kellner}
\date{}
\maketitle
\vspace*{-0.8cm}

\abstract{Let $B_n$ ($n = 0, 1, 2, \ldots$) denote the usual $n$-th Bernoulli number.
Let $l$ be a positive even integer where $l=12$ or $l \geq 16$.
It is well known that the numerator of the reduced quotient $|B_l/l|$ is
a product of powers of irregular primes.
Let $(p,l)$ be an irregular pair with $B_l/l \not\equiv B_{l+p-1}/(l+p-1) \modp{p^2}$.
We show that for every $r \geq 1$ the congruence $B_{m_r}/m_r \equiv 0 \modp{p^r}$
has a unique solution $m_r$ where $m_r \equiv l \modp{p-1}$ and
$l \leq m_r < (p-1)p^{r-1}$. The sequence $(m_r)_{r \geq 1}$
defines a $p$-adic integer $\chi_{(p,\,l)}$ which is a zero
of a certain $p$-adic zeta function $\zeta_{p,\,l}$ originally defined by
T.~Kubota and H.~W.~Leopoldt. We show some properties of these functions
and give some applications.
Subsequently we give several computations of the
(truncated) $p$-adic expansion of $\chi_{(p,\,l)}$ for irregular pairs $(p,l)$
with $p$ below 1000.}

\smallskip

\vspace*{-0.5cm}

\section{Introduction}

The classical Bernoulli numbers $B_n$ are defined by
\[
  \frac{z}{e^z-1} = \sum_{n=0}^\infty B_n \frac{z^n}{n!} ,
    \quad |z| < 2 \pi .
\]
These numbers are rational where $B_n = 0$ for odd $n > 1$ and
$(-1)^{\frac{n}{2}+1}B_n > 0$ for even $n > 0$.
The Bernoulli numbers are connected with the Riemann zeta function $\zeta(s)$
on the positive real axis by Euler's formula
\begin{equation} \label{eqn-zeta-n}
   \zeta( n ) = -\frac{1}{2} \frac{(2 \pi i)^n}{n!} \, B_n ,
     \quad n \in \NN , \,\, 2 \pdiv n ;
\end{equation}
the functional equation of $\zeta(s)$ leads to
\begin{equation} \label{eqn-zeta-1-n}
   \zeta( 1-n ) = - \frac{ B_{n} }{n} \quad \mbox{for} \quad n \in \NN , \,\, n \geq 2 .
\end{equation}

Let, as usual, $\eulerphi$ denote Euler's totient function.
The Kummer congruences state
for $n, m, p, r \in \NN$, $n, m$ even, $p$ prime with $p-1 \notdiv n$ that
\begin{equation} \label{eqn-kummer-congr}
   (1-p^{n-1}) \frac{B_n}{n} \equiv (1-p^{m-1}) \frac{B_{m}}{m} \pmod{p^r}
\end{equation}
when $n \equiv m \modp{\eulerphi(p^r)}$; see \cite[Thm.~5, p.~239]{ireland90}.
\pagebreak

In 1850 E.\,E.~Kummer \cite{kummer50} introduced two classes of odd primes,
later called regular and irregular.
An odd prime $p$ is called \textsl{regular} if $p$ does not divide the class number of
the cyclotomic field $\QQ(\mu_p)$ where $\mu_p$ is the set of $p$-th roots of unity;
other\-wise \textsl{irregular}. Kummer proved that Fermat's last theorem (FLT) is true
if the exponent is a regular prime. Kummer also gave an equivalent definition of irregularity concerning
Bernoulli numbers. We recall the usual definition from \cite[pp.~233--234]{ireland90}.

\begin{defin} \label{def-irrpair}
Let $p$ be an odd prime.
The pair $(p,l)$ is called an \textsl{irregular pair} if $p$ divides the numerator
of $B_l$ where $l$ is even and $2 \leq l \leq p-3$.
The \textsl{index of irregularity} of $p$ is defined to be
\[
    i(p) := \# \{ (p,l) \mbox{\ is an irregular pair} :\, l = 2, 4, \ldots, p-3 \} .
\]
The prime $p$ is called a \textsl{regular} prime if $i(p) = 0$, otherwise
an \textsl{irregular} prime.
\end{defin}

We introduce the following notations for rational numbers. If $q$ is rational then we use the
representation $q=N/D$ where $(N,D)=1$ and $D > 0$. We define $\denom(q) = D$ resp.\
$\numer(q) = N$ for the denominator resp.\ the numerator of $q$. The notation $m \pdiv q$
where $m$ is a positive integer means that $m \pdiv \numer(q)$; we shall also write
$q \equiv 0 \modp{m}$ in this case.

For now, let $n$ be an even positive integer.
An elementary property of the Ber\-noulli numbers, independently discovered
by T.~Clausen \cite{clausen40} and K.\,G.\,C.\ von Staudt \cite{staudt40} in 1840, is the following.
The structure of the denominator of $B_n$ is given by
\begin{equation} \label{eqn-clausen-von-staudt}
   B_n + \sum_{p-1 \pdiv n} \frac{1}{p} \in \ZZ \qquad \mbox{which implies} \qquad
     \denom(B_n) = \prod_{p-1 \pdiv n} p .
\end{equation}

A further result, often associated with the name of
J.\,C.~Adams, see \cite[Prop.\ 15.2.4, p.~238]{ireland90} and Section \ref{sect-adam} below,
is that $B_n/n$ is a $p$-integer for all primes $p$ with $p-1 \notdiv n$. Therefore
\begin{equation} \label{eqn-adams-triv-div}
   \prod_{p-1 \notdiv n} p^{\ord_p n}
\end{equation}
divides $\numer(B_n)$; since this factor is cancelled in the numerator of $B_n/n$,
we shall call it the \textsl{trivial factor} of $B_n$.
By $|B_n| > 2(n \big/ {2\pi e})^n$, see \cite[Eq.\ (8), p.~232]{ireland90},
and Table \ref{tab-bn}, one can easily show that the numerator of
$|B_n/n|$ equals 1 for $2 \leq n \leq 10$ and $n=14$.
Otherwise this numerator is a product of powers of irregular primes;
this is a consequence of the Kummer congruences.
The determination of irregular primes resp.\ irregular pairs is
still a difficult task, see \cite{irrprime12M}.
One can easily show that infinitely many irregular primes exist;
for a short proof see Carlitz \cite{carlitz54}.
Unfortunately, the more difficult question of
whether infinitely many regular primes exist is still open.
However, calculations in \cite{irrprime12M} show that about 60\% of all primes less than
12 million are regular which agree with an expected distribution proposed by
Siegel \cite{siegel64}.

The values of $B_n$ and $B_n/n$ for $n \leq 20$ are given in Table \ref{tab-bn},
irregular pairs with $p < 1000$ in Table \ref{tab-irrpair-ord-10}.
For brevity we write $\BN(n) = B_n / n$; these are called the
\textsl{divided} Bernoulli numbers. Throughout this paper
all indices concerning Bernoulli numbers are even and $p$ denotes an odd prime.

\subsection*{Acknowledgement}
The author wishes to thank Prof.\ Samuel J.\ Patterson
for his numerous helpful comments and suggestions to improve this paper.
The author is thankful to Prof.\ L.\,C.\ Washington for some comments on Iwasawa theory
and reference \cite{pollaczek24}. The author is also grateful to the referee for useful
suggestions and detailed remarks. This paper is based on the author's diploma thesis
\cite{kellner02}.

\section{Preliminaries}
\setcounter{equation}{0}

The definition of irregular pairs can be extended to irregular prime powers
as was first proposed by the author \cite[Section 2.5]{kellner02}.

\begin{defin} \label{def-ordn-irrpair}
Let $p$ be an odd prime and $n$, $l$ be positive integers.
A pair $(p,l)$ is called an \textsl{irregular pair of order $n$} if
$p^n \pdiv \BN(l)$ where $l$ is even and $2 \leq l < \eulerphi(p^n)$.
Define
\[
   \IRR_n := \{ (p, l) :\, p^n \pdiv \BN(l), \, p \mbox{\ is an odd prime}, \,
     2 \leq l < \eulerphi(p^n), \, 2 \pdiv l \}
\]
as the set of irregular pairs of order $n$. For a prime $p$ the \textsl{index}
of irregular pairs of order $n$ is defined by
\[
    i_n(p) := \# \{ (p,l) :\, (p,l) \in \IRR_n \} .
\]
Define the map
\[
    \lambda_n: \,\, \IRR_{n+1} \to \IRR_{n} ,
      \quad (p,l) \mapsto (p,l \modres \eulerphi(p^n))
\]
where $x \modres y$ gives the least nonnegative residue $x'$ with $0 \leq x' < y$
for positive integers $x$ and $y$.
Let $(p,l_n) \in \IRR_n$ and $(p,l_m) \in \IRR_m$
be irregular pairs of order $n$ resp.\ $m$ where $n > m \geq 1$.
We say that $(p,l_n)$ is \textsl{related} to $(p,l_m)$ if
$l_n \equiv l_m \modp{\eulerphi(p^m)}$ holds.
Note that ``related'' is not a symmetric relation.
\end{defin}

\begin{remark} \label{rem-ordn-irrpair}
This definition includes the older
Definition \ref{def-irrpair} in the case $n=1$.
Therefore one has $i_1(p) = i(p)$. Let $(p,l) \in \IRR_n$
with $n \geq 1$.
The Kummer congruences \refeqn{eqn-kummer-congr} imply
that $p^n \pdiv \BN(l+\nu \eulerphi(p^n))$ for all $\nu \in \NN_0$.
Conversely, if $p^n \pdiv \BN(m)$ for some even integer $m$, then
there exists an irregular pair $(p,l) \in \IRR_n$
where $l \equiv m \modp{\eulerphi(p^n)}$ with $l \leq m$ holds.

The map $\lambda_n$ is well defined by the properties mentioned above.
Let $(p,l_n) \in \IRR_n$ and $(p,l_m) \in \IRR_m$ where $n > m \geq 1$
and $(p,l_n)$ is related to $(p,l_m)$.
By applying the maps $\lambda_{n-1}, \lambda_{n-2}, \ldots, \lambda_{m}$
one derives a chain of related irregular pairs of descending order
\begin{equation} \label{eqn-chain-irrpair}
   (p,l_n) \in \IRR_n, \ (p,l_{n-1}) \in \IRR_{n-1}, \
     (p,l_{n-2}) \in \IRR_{n-2}, \ \ldots, \ (p,l_m) \in \IRR_m
\end{equation}
where
\[
    l_n \geq l_{n-1} \geq l_{n-2} \geq \ldots \geq l_m .
\]
\end{remark}

\begin{defin} \label{def-delta}
For $(p,l) \in \IRR_n$, $n \geq 1$ define
\[
    \Delta_{(p,\,l)} \, \equiv \, p^{-n} \left(
       \BN(l + \eulerphi(p^n)) - \BN(l) \right) \pmod{p}
\]
with $0 \leq \Delta_{(p,\,l)} < p$.
When $\Delta_{(p,\,l)} = 0$ we call $\Delta_{(p,\,l)}$
\textsl{singular}. For an irregular prime $p$ set
\[
   \Delta(p) := \left\{
     \begin{array}{rl}
       1 , & \Delta_p \neq 0 \\
       0 , & \Delta_p = 0
     \end{array} \right.
\]
with
\[
    \Delta_p = \prod_{\nu=1}^{i(p)} \Delta_{(p,l_\nu)} ,
      \quad  (p,l_\nu) \in \IRR_1 .
\]
Then $\Delta(p) = 1$ if and only if all $\Delta_{(p,l_\nu)}$ are nonsingular.
\end{defin}

We need a generalized form of the Kummer congruences
which allows us to obtain most of the later results,
see Carlitz \cite[Thm.~3, p.~425]{carlitz53},
especially Fresnel \cite[Cor.~6, p.~319]{fresnel67}.

\begin{theorem}[Carlitz] \label{theor-kummer-congr-general}
Let $k, m, n, p, r, \omega \in \NN$, $m$ even, $p$ prime with $p-1 \notdiv m$,
and $\omega = k \eulerphi(p^n)$. Then
\begin{equation} \label{eqn-kummer-congr-general}
   \sum_{\nu=0}^r \binom{r}{\nu} (-1)^{\nu}
     (1-p^{m+ \nu \omega -1}) \BN(m+\nu \omega) \equiv 0 \pmod{p^{nr}} .
\end{equation}
\end{theorem}

Here we need a special version without Euler factors which $p$-adically
shows the periodic behavior of the divided Bernoulli numbers.

\begin{corl} \label{corl-kummer-congr-irrpair}
Let $(p,l) \in \IRR_n$, $n \geq 1$.
Let $k, m, r, \omega \in \NN$, $r > 1$, and $\omega = k \eulerphi(p^n)$.
For $m=l+j \eulerphi(p^n)$ with $j \geq 0$ we have
\[
   \sum_{\nu=0}^r \binom{r}{\nu} (-1)^{\nu} \,
     p^{-n} \BN(m+\nu \omega) \equiv 0
      \pmod{(p^{m-1},p^{n(r-1)})} .
\]
\end{corl}

\begin{proof}
Since $(p,l) \in \IRR_n$ we know that $p^n \pdiv \BN(l+j\eulerphi(p^n))$ for all
$j \geq 0$. Hence, we can reduce Congruence \refeqn{eqn-kummer-congr-general}
to $\mods{p^{n(r-1)}}$ when multiplying it by $p^{-n}$.
One easily sees that all Euler factors in the sum of \refeqn{eqn-kummer-congr-general}
vanish $\mods{p^{m-1}}$.
\end{proof}

Proposition \ref{prop-irrpair-n} below shows how to find irregular pairs of
higher order. Beginning from an irregular pair of order $n$, we can
characterize related irregular pairs of order $n+1$ if they exist.
First we need a lemma.

\begin{lemma} \label{lemma-seq-a2}
Let $n$ be a positive integer and $p$ be an odd prime.
Let $( \alpha_\nu )_{\nu \geq 0}$ be a sequence of $p$-integers $\alpha_\nu \in \QQ$
for all $\nu \in \NN_0$. If one has
\begin{equation} \label{eqn-seq-a2-mod}
   \alpha_\nu - 2 \alpha_{\nu+1} + \alpha_{\nu+2} \equiv 0 \pmod{p^n}
\end{equation}
then the sequence is equidistant $\mods{p^n}$. For
$\alpha_0 \not\equiv \alpha_1 \modp{p}$ the elements $\alpha_0$ up to $\alpha_{p^n-1}$
run through all residues $\mods{p^n}$. Then exactly one element
$\alpha_s \equiv 0 \modp{p^n}$ exists with $0 \leq s < p^n$
where $s \equiv -\alpha_0(\alpha_1-\alpha_0)^{-1} \modp{p^n}$.
\end{lemma}

\begin{proof}
We rewrite Congruence \refeqn{eqn-seq-a2-mod} to
\[
   \alpha_\nu - \alpha_{\nu+1} \equiv \alpha_{\nu+1} - \alpha_{\nu+2} \pmod{p^n}
     , \quad \nu \in \NN_0 .
\]
One easily observes by induction on $\nu$ that all elements $\alpha_\nu$ are
equidistant $\mods{p^n}$. Let $\delta \equiv \alpha_1 - \alpha_0 \modp{p^n}$,
then we obtain
\[
   \alpha_\nu \equiv \alpha_0 + \delta \nu \pmod{p^n} .
\]
The case $\alpha_0 \not\equiv \alpha_1 \modp{p}$
provides that $\alpha_0 + \delta \nu$ resp.\ $\alpha_\nu$ runs
through all residues $\mods{p^n}$ for
$0 \leq \nu < p^n$, since $\delta$ is invertible $\mods{p^n}$.
Then exactly one element $\alpha_s$ exists with
$0 \equiv \alpha_s \equiv \alpha_0 + \delta s \modp{p^n}$ and $0 \leq s < p^n$.
This finally gives
$s \equiv -\alpha_0 \, \delta^{-1} \equiv -\alpha_0(\alpha_1-\alpha_0)^{-1} \modp{p^n}$.
\end{proof}

\begin{prop} \label{prop-irrpair-n}
Let $(p,l) \in \IRR_n$, $n \geq 1$.
Let the sequence $( \alpha_j )_{j \geq 0}$ satisfy
\[
    \alpha_j \equiv p^{-n} \BN(l + j \eulerphi(p^n))  \pmod{p} .
\]
For $\Delta_{(p,\,l)} \equiv \alpha_1 - \alpha_0 \modp{p}$
where $0 \leq \Delta_{(p,\,l)} < p$ there exist three cases:
\begin{enumerate}
 \item If $\Delta_{(p,\,l)} = 0$ and $\alpha_0 \not\equiv 0 \modp{p}$,
   then there are no related irregular pairs of order $n+1$ and higher,
 \item If $\Delta_{(p,\,l)} = 0$ and $\alpha_0 \equiv 0 \modp{p}$,
   then all $(p,l+\nu \eulerphi(p^n)) \in \IRR_{n+1}$
   are related irregular pairs of order $n+1$ for $\nu=0,\ldots,p-1$,
 \item If $\Delta_{(p,\,l)} \neq 0$, then exactly one related
   irregular pair of order $n+1$ exists. One has
   $(p,l+s \eulerphi(p^n)) \in \IRR_{n+1}$ with $0 \leq s < p$ where
   $s \equiv -\alpha_0 \Delta_{(p,\,l)}^{-1} \modp{p}$.
\end{enumerate}
\end{prop}

\begin{proof}
Let $j \geq 0$.
Using Corollary \ref{corl-kummer-congr-irrpair}
with $r=2$, $\omega = \eulerphi(p^n)$, and $m=l+j \eulerphi(p^n) \geq 2$,
we get
\[
   \sum_{\nu=0}^2 \binom{2}{\nu} (-1)^{\nu} \,
     p^{-n} \BN(m+\nu \omega) \equiv 0
      \pmod{p}
\]
which is
\[
   \alpha_j - 2 \alpha_{j+1} + \alpha_{j+2} \equiv 0 \pmod{p} .
\]
This satisfies the conditions of Lemma \ref{lemma-seq-a2}. We obtain
three cases: \smallskip

Case (1): We have $\alpha_0 \equiv \alpha_1 \modp{p}$ and
$\alpha_0 \not\equiv 0 \modp{p}$. One observes that
$\alpha_j \not\equiv 0 \modp{p}$ resp.\ $p^{n+1} \notdiv \BN(l + j \eulerphi(p^n))$
for all $j \geq 0$.
Therefore, there cannot exist related irregular pairs of order $n+1$.
Also there cannot exist related irregular pairs of order $r > n+1$,
otherwise we would get a contradiction to \refeqn{eqn-chain-irrpair}.
\smallskip

Case (2): We have $\alpha_0 \equiv \alpha_1 \modp{p}$ and
$\alpha_0 \equiv 0 \modp{p}$. This yields
$\alpha_j \equiv 0 \modp{p}$ resp.\ $p^{n+1} \pdiv
\BN(l + j \eulerphi(p^n))$ for all $j \geq 0$.
Hence, $p$ related irregular pairs of order $n+1$ exist where
$(p,l+\nu \eulerphi(p^n)) \in \IRR_{n+1}$ for $\nu=0,\ldots,p-1$.
\smallskip

Case (3): We have $\alpha_0 \not\equiv \alpha_1 \modp{p}$.
Lemma \ref{lemma-seq-a2} provides exactly one element
$\alpha_s \equiv 0 \modp{p}$ with the desired properties.
Hence, $(p,l+s \eulerphi(p^n))$ is the only related irregular pair
of order $n+1$.
\end{proof}

\begin{remark}
Vandiver \cite{vandiver37} describes the result of the previous
proposition for the case $n=1$ and only for the first irregular primes
37, 59, and 67. For these primes Pollaczek \cite{pollaczek24} has calculated
the indices $s$ of the now called irregular pair of order two,
but case $p=67$ with $s=2$ is incorrect, see column $s_2$ of Table \ref{tab-irrpair-ord-10}.
This error was already noticed by Johnson \cite{johnson74} who has also
determined all irregular pairs $(p,l')$ of order two with $p$ below 8000.
Wagstaff \cite{wagstaff78} has extended calculations of irregular pairs,
indices $s$, and associated cyclotomic invariants up to $p < 125\,000$.
He also checked that FLT is true for all such exponents $p$ in that range.
Finally, Buhler, Crandall, Ernvall, Mets\"ankyl\"a, and Shokrollahi \cite{irrprime12M}
have extended calculations of irregular pairs and associated
cyclotomic invariants up to $p < 12\,000\,000$.
For all these irregular pairs $(p,l)$ in that range
$\Delta_{(p,\,l)} \neq 0$ is always valid
which ensures that each time there is only one related irregular pair $(p,l')$
of order two. Hence $i_2(p) = i(p)$ for these irregular primes $p$.
One has to notice that always $(p,l) \neq (p,l')$. So far, no irregular pair $(p,l)$
has been found with $p^2 \pdiv \BN(l)$.
\end{remark}

Using Proposition \ref{prop-irrpair-n} one can successively find irregular pairs
of higher order. We can easily extend the result starting from an irregular pair
$(p,l) \in \IRR_n$ and requiring that $l > n$
to obtain a related irregular pair $(p,l') \in \IRR_{2n}$.

\begin{prop} \label{prop-irrpair-2n}
Let $(p,l) \in \IRR_n$, $n \geq 1$.
Suppose that $l > n$. Let the sequence $( \alpha_j )_{j \geq 0}$ satisfy
\[
    \alpha_j \equiv p^{-n} \BN(l + j \eulerphi(p^n))  \pmod{p^n} .
\]
If $\Delta_{(p,\,l)} \neq 0$, then there is exactly one related irregular pair of order $2n$
\[
   (p,l+s \eulerphi(p^n)) \in \IRR_{2n}
\]
with $0 \leq s < p^n$ where
$s \equiv -\alpha_0 (\alpha_1-\alpha_0)^{-1} \modp{p^n}$.
Correspondingly, there also exists exactly one related irregular pair of order $n+\nu$
\[
   (p,l+s_\nu \eulerphi(p^n)) \in \IRR_{n+\nu}
\]
for each $\nu = 1,\ldots,n-1$ with
$0 \leq s_\nu < p^\nu$ where $s_\nu \equiv s \modp{p^\nu}$.
\end{prop}

\begin{proof}
Let $j \geq 0$.
Using Corollary \ref{corl-kummer-congr-irrpair} again
with $r=2$, $\omega = \eulerphi(p^n)$, and $m=l+j \eulerphi(p^n) > n$
yields
\[
   \alpha_j - 2 \alpha_{j+1} + \alpha_{j+2} \equiv 0 \pmod{p^n} .
\]
If $\Delta_{(p,\,l)} \neq 0$, then Lemma \ref{lemma-seq-a2}
provides exactly one element $\alpha_s \equiv 0 \modp{p^n}$
with $0 \leq s < p^n$ where $s \equiv -\alpha_0(\alpha_1-\alpha_0)^{-1} \modp{p^n}$.
Therefore, $(p,l+s \eulerphi(p^n))$ is the only related irregular pair of order $2n$.
Similarly, regarding the congruences above $\mods{p^\nu}$ instead of $\mods{p^n}$
for $\nu=1,\ldots,n-1$ yields the proposed properties.
\end{proof}

Finally, we start from an irregular pair $(p,l) \in \IRR_n$
where we have to suppose that $l > (r-1)n$ with $r \geq 2$
to obtain a related irregular pair $(p,l') \in \IRR_{rn}$.

\begin{prop} \label{prop-irrpair-rn}
Let $(p,l) \in \IRR_n$, $n \geq 1$.
Let $r>1$ be an integer and suppose that $l>(r-1)n$.
Let the sequence $( \alpha_j )_{j \geq 0}$ satisfy
\[
    \alpha_j \equiv p^{-n} \BN(l + j \eulerphi(p^n))  \pmod{p^{(r-1)n}} .
\]
Then this sequence satisfies for all $j \geq 0$:
\[
   \sum_{\nu=0}^r \binom{r}{\nu} (-1)^{\nu}
      \alpha_{\nu+j} \equiv 0
      \pmod{p^{(r-1)n}} .
\]
The elements $\alpha_0$ up to $\alpha_{r-1}$ induce the entire
sequence $( \alpha_j )_{j \geq 0}$. Elements with
$\alpha_s \equiv 0 \modp{p^{(r-1)n}}$ where $0 \leq s < p^{(r-1)n}$
provide related irregular pairs of order $rn$ with $(p,l+s \eulerphi(p^n)) \in \IRR_{rn}$.
If $\Delta_{(p,\,l)} \neq 0$ and the elements
$\alpha_0$ up to $\alpha_{r-1}$ are equidistant $\mods{p^{(r-1)n}}$,
then there is exactly one related irregular pair of order $rn$
\[
   (p,l+s \eulerphi(p^n)) \in \IRR_{rn}
\]
with $0 \leq s < p^{(r-1)n}$ where
$s \equiv -\alpha_0 (\alpha_1-\alpha_0)^{-1} \modp{p^{(r-1)n}}$.
Correspondingly, there exists exactly one related irregular pair of order $n+k$
\[
   (p,l+s_k \eulerphi(p^n)) \in \IRR_{n+k}
\]
for each $k = 1,\ldots,(r-1)n-1$
with $0 \leq s_k < p^k$ where $s_k \equiv s \modp{p^k}$.
\end{prop}

\begin{proof}
Let $j \geq 0$.
Clearly, by Corollary \ref{corl-kummer-congr-irrpair}
with $r>1$, $\omega = \eulerphi(p^n)$, and $m=l+j \eulerphi(p^n) > (r-1)n$,
we have
\[
   \sum_{\nu=0}^r \binom{r}{\nu} (-1)^{\nu}
      \alpha_{\nu+j} \equiv 0  \pmod{p^{(r-1)n}} .
\]
This induces the whole sequence $( \alpha_j )_{j \geq 0}$ by
\begin{equation} \label{eqn-loc-irrpair-rn-1}
   \alpha_{r+j} \equiv (-1)^{r+1} \sum_{\nu=0}^{r-1} \binom{r}{\nu} (-1)^{\nu}
      \alpha_{\nu+j} \pmod{p^{(r-1)n}} .
\end{equation}
Among all elements $\alpha_s$ with $0 \leq s < p^{(r-1)n}$,
an element $\alpha_s \equiv 0 \modp{p^{(r-1)n}}$ provides
a related irregular pair $(p,l+s \eulerphi(p^n)) \in \IRR_{rn}$ of order $rn$.
\smallskip

Now we assume that $\Delta_{(p,\,l)} \neq 0$ and the first elements
$\alpha_0$ up to $\alpha_{r-1}$ are equidistant $\mods{p^{(r-1)n}}$.
We show that this property transfers to the entire sequence.
Let $\gamma, \delta$ be integers. It is easily seen for $r>1$ that
\begin{equation} \label{eqn-loc-irrpair-rn-2}
   \sum_{\nu=0}^r \binom{r}{\nu} (-1)^{\nu}
      (\gamma + \delta \nu ) = 0 .
\end{equation}
Consider $\gamma \equiv \alpha_0 \modp{p^{(r-1)n}}$
and $\delta \equiv \alpha_1-\alpha_0 \modp{p^{(r-1)n}}$ where $p \notdiv \delta$ by assumption.
Combining \refeqn{eqn-loc-irrpair-rn-1} and \refeqn{eqn-loc-irrpair-rn-2} yields
\begin{equation} \label{eqn-loc-irrpair-rn-3}
   (-1)^{r+1} \sum_{\nu = 0}^{r-1} \binom{r}{\nu} (-1)^{\nu} \, (\alpha_0 + \delta \nu)
     \equiv \alpha_0 + \delta r \equiv \alpha_r  \pmod{p^{(r-1)n}}
\end{equation}
which shows inductively that all successive elements $\alpha_j$ with $j \geq r$
are equidistant $\mods{p^{(r-1)n}}$.
Since $\delta$ is invertible $\mods{p^{(r-1)n}}$, exactly one solution exists
with $0 \equiv \alpha_s \equiv \alpha_0 + \delta s \modp{p^{(r-1)n}}$ where
$s \equiv -\alpha_0 (\alpha_1-\alpha_0)^{-1} \modp{p^{(r-1)n}}$.
Using similar arguments, Congruences \refeqn{eqn-loc-irrpair-rn-1} and
\refeqn{eqn-loc-irrpair-rn-3} are also valid $\mods{p^k}$ for $k=1,\ldots,(r-1)n-1$
which provides for each $k$ a unique solution $s_k$ with $0 \leq s_k < p^k$ where
$s_k \equiv s \modp{p^k}$ holds.
\end{proof}

In \cite[pp.~125--130]{kellner02} several examples and calculations are given
which use the previous propositions. These results are reprinted in
Table \ref{tab-calc-37-1} and following.
It will turn out later that calculations can be
further simplified. Because of the rapid growth of indices,
it is useful to write indices of irregular pairs of
higher order $p$-adically.

\begin{defin} \label{def-pair-padic}
Let $(p,l) \in \IRR_n$, $n \geq 1$. We write
\[
   (p,s_1,s_2, \ldots, s_n) \in \IRP_n \quad \mbox{where} \quad
     l = \sum_{\nu=1}^{n} s_\nu \, \eulerphi( p^{\nu-1} )
\]
for the $p$-adic notation of $(p,l)$
with $0 \leq s_\nu < p$ for $\nu = 1,\ldots,n$ and $2 \leq s_1 \leq p-3$, $2 \pdiv s_1$.
The corresponding set is denoted as $\IRP_n$,
the map corresponding to $\lambda_n$ is given by
\[
    \widehat{\lambda}_n: \,\, \IRP_{n+1} \to \IRP_{n},
      \quad (p,s_1,s_2, \ldots, s_n,s_{n+1}) \mapsto (p,s_1,s_2, \ldots, s_n) .
\]
The pair $(p,l)$ and the element $(p,s_1,s_2, \ldots, s_n)$ are called \textsl{associated}.
\end{defin}

\begin{remark}
The definition of $\IRP_n$ means that we have $\IRR_1 = \IRP_1$ for $n=1$.
For $n \geq 2$ we can define a map $\IRR_n \to \IRP_n$,
$(p,l) \mapsto (p,s_1,\ldots,s_n)$ where the $s_k$
are uniquely determined by the $p$-adic representation
\[
   l = s_1 + (p-1) \hat{s} , \quad \hat{s} =
     \sum_{\nu=0}^{n-2} s_{\nu+2} \, p^\nu,
       \quad 0 \leq s_{\nu+2} < p
\]
and by $s_1 \equiv l \modp{p-1}$ with $2 \leq s_1 \leq p-3$.
If $s_k = 0$ with $k \geq 2$ then there is
an irregular pair $(p,l_k)$ of order $k$ with
\[
    (p,l_k) \in \IRR_k \quad \mbox{and} \quad  (p,l_k) \in \IRR_{k-1} .
\]
Note that $(p,s_1,s_2, \ldots, s_n)$ is also called an \textsl{irregular pair}
with $(s_1,s_2, \ldots, s_n)$ as the second parameter given $p$-adically.
\end{remark}

\section{Main results}
\setcounter{equation}{0}
\label{sect-main}

\begin{theorem} \label{theor-delta-irr}
Let $(p,l_1) \in \IRR_1$. If $\Delta_{(p,l_1)} \neq 0$ then
for each \mbox{$n > 1$} there exists exactly one related irregular pair of order $n$.
There is a unique sequence $(l_n)_{n \geq 1}$ resp.\ $(s_n)_{n \geq 1}$ with
\[
   (p,l_n) \in \IRR_n \quad \mbox{resp.} \quad (p,s_1,\ldots,s_n) \in \IRP_n
\]
and
\[
    l_1 \leq l_2 \leq l_3 \leq \ldots \,,
      \quad \lim_{n \to \infty} l_n = \infty .
\]
Moreover one has
\[
   \Delta_{(p,l_1)} = \Delta_{(p,l_2)} = \Delta_{(p,l_3)} = \ldots \,.
\]
If $\Delta(p)=1$ then
\[
   i(p) = i_2(p) = i_3(p) = \ldots \,.
\]
\end{theorem}

\begin{theorem} \label{theor-delta-irr-0}
Let $(p,l_n) \in \IRR_n$, $n \geq 1$
with $\Delta_{(p,l_n)}=0$. Then there are two cases:
\begin{enumerate}
 \item $(p,l_n) \notin \IRR_{n+1}$$:$
   There are no related irregular pairs of order $n+1$ and higher,
 \item $(p,l_n) \in \IRR_{n+1}$$:$
   There exist $p$ related irregular pairs of order $n+1$ where
   $(p,l_{n+1,j}) \in \IRR_{n+1}$
   with $\Delta_{(p, \,l_{n+1,j})} = 0$ and $l_{n+1,j} = l_n+j\eulerphi(p^n)$
   for $j=0,\ldots,p-1$.
\end{enumerate}
\end{theorem}

The property of $\Delta_{(p,\,l)}$, whether $\Delta_{(p,\,l)}$ vanishes or not,
is passed on to all related irregular pairs of higher order.
The case of a singular $\Delta_{(p,\,l)}$ would possibly imply
a strange behavior without any regularity.
By calculation in \cite{irrprime12M} up to $p < 12\,000\,000$, no
such $\Delta_{(p,\,l)}$ was found. The following diagram illustrates
both cases.

\begin{diagenv} \label{diag-tree}
\end{diagenv}

\unitlength0.0075cm
\begin{center}
\begin{picture}(700,520)
  \put(200,460){$\Delta_{(p,\,l_1)} = 0$}
  \put(10,390){$\IRR_1$}
  \put(10,290){$\IRR_2$}
  \put(10,190){$\IRR_3$}
  \put(10,090){$\IRR_4$}
  \put(145,030){$l_1$}
  \put(245,030){$l_2$}
  \put(395,030){$l_3$}
  \put(150,400){\line(0,-1){100}}
  \put(150,400){\line(1,-1){100}}
  \put(150,400){\line(2,-1){200}}
  \put(150,400){\circle*{10}}
  \put(150,300){\circle*{10}}
  \put(250,300){\circle*{10}}
  \put(350,300){\circle*{10}}

  \put(250,300){\line(0,-1){100}}
  \put(250,300){\line(3,-4){75}}
  \put(250,300){\line(3,-2){150}}

  \put(325,200){\circle*{10}}
  \put(250,200){\circle*{10}}
  \put(400,200){\circle*{10}}

  \put(250,200){\line(0,-1){100}}
  \put(250,200){\line(1,-2){50}}
  \put(250,200){\line(1,-1){100}}
  \put(400,200){\line(0,-1){100}}
  \put(400,200){\line(1,-2){50}}
  \put(400,200){\line(1,-1){100}}

  \put(250,100){\circle*{10}}
  \put(300,100){\circle*{10}}
  \put(350,100){\circle*{10}}
  \put(400,100){\circle*{10}}
  \put(450,100){\circle*{10}}
  \put(500,100){\circle*{10}}
\end{picture}
\begin{picture}(500,520)
  \put(180,460){$\Delta_{(p,\,l_1)} \neq 0$}
  \put(10,390){$\IRR_1$}
  \put(10,290){$\IRR_2$}
  \put(10,190){$\IRR_3$}
  \put(10,090){$\IRR_4$}
  \put(145,030){$l_1$}
  \put(220,030){$l_2,\,l_3$}
  \put(345,030){$l_4$}
  \put(150,400){\line(1,-1){100}}
  \put(150,400){\circle*{10}}
  \put(250,300){\circle*{10}}
  \put(250,300){\line(0,-1){100}}
  \put(250,200){\circle*{10}}
  \put(250,200){\line(1,-1){100}}
  \put(350,100){\circle*{10}}
\end{picture}
\end{center}

Here a vertical line indicates that $(p,l_n) \in \IRR_n \cap \IRR_{n+1}$ happens.
On the left side we then have $p$ related irregular pairs of order $n+1$ which are
represented by branches. In this case the corresponding Bernoulli number $\BN(l_n)$
\textsl{decides} whether there exist further branches or they stop.
Instead of $n$ the order of the $p$-power must be at least $n+1$.
This also means that an associated irregular pair $(p,s_1,\ldots,s_{n+1}) \in \IRP_{n+1}$ must
have a zero $s_{n+1}=0$ each time.
In contrast to, the right side shows that in the case of $\Delta_{(p,l_1)} \neq 0$
there is only one related irregular pair of each higher order.
If $\Delta(p)=1$ then higher powers $p^\nu$ are equally distributed
among the numerators of $\BN(n)$. For each irregular pair considered,
there exists exactly one index $n_{k,\nu}$ with $n_{k,\nu} = n_{0,\nu} + k\eulerphi(p^\nu)$,
$k \in \NN_0$ in the disjoint intervals
\[
    (k \, \eulerphi(p^\nu), (k+1) \, \eulerphi(p^\nu) )
\]
for which $p^\nu \pdiv \BN(n_{k,\nu})$ is valid.
\smallskip

In \cite[pp.~128--130]{kellner02} irregular pairs of order 10 were calculated
for all irregular primes $p < 1000$. These results are
reprinted in Table \ref{tab-irrpair-ord-10}. In this table
only one irregular pair has a zero in its $p$-adic notation:
\[
    (157,62,40,145,67,29,69,0,87,89,21) \in \IRP_{10} .
\]
Hence, one has with a \textsl{relatively small} index that
\[
    (157, 6\,557\,686\,520\,486) \in \IRR_6 \cap \IRR_7 .
\]
It seems that these zeros can be viewed as exceptional;
see also Table \ref{tab-irrpair-ord-100}.
It would be of interest to investigate in which regions such indices may occur.
This could explain why no irregular pair $(p,l) \in \IRR_1 \cap \IRR_2$
has yet been found, because these regions are beyond present calculations.
Here we have index $12\,000\,000$ in \cite{irrprime12M} against
index $6\,557\,686\,520\,486$.
Because of the rare occurrence of zeros one can expect
that $(p,l) \in \IRR_1 \cap \IRR_2$ resp.\ $p^2 \pdiv \BN(l)$ will
not happen often.
\smallskip

Another phenomenon is the occurrence of huge irregular prime factors.
Wagstaff \cite{wagstaff00} has completely factored
the numerators of the Bernoulli numbers $B_n$ with index up to
$n=152$. Most of these irregular prime factors are
large numbers, the greatest factors have 30 up to 100 digits.
\smallskip

Finally, summarizing all facts together,
the property $\Delta(p)=1$ can be viewed as a structural property of the Bernoulli numbers.
This leads us to the following conjecture named by the author $\Delta$-Conjecture.

\begin{conj}[$\Delta$-Conjecture] \label{conj-delta}
For all irregular primes $p$ the following properties,
which are equivalent, hold:
\begin{enumerate}
\item $\Delta_{(p,\,l)}$ is nonsingular for all irregular pairs $(p,l) \in \IRR_1$,
\item $\Delta(p)=1$,
\item $i(p)=i_2(p)=i_3(p)=\ldots\,$.
\end{enumerate}
\end{conj}

Assuming the $\Delta$-Conjecture one can also prove the existence of infinitely many
irregular primes using only information about the numerators of divided Bernoulli
numbers, see \cite[Satz 2.8.2, p.~87]{kellner02}.
Now we give the proofs of the theorems above.

\begin{prop} \label{prop-irr-delta-n}
Let $(p,l_n) \in \IRR_n$, $n \geq 1$ with
$\Delta_{(p,l_n)} \neq 0$. Then there is exactly one related
irregular pair $(p,l_{n+1}) \in \IRR_{n+1}$ with
$\Delta_{(p,l_n)} = \Delta_{(p,l_{n+1})}$.
\end{prop}

\begin{proof}
We write $\Delta_n = \Delta_{(p,l_n)}$. Note that $l_n > 2$ and $p > 3$. Define the sequence
$( \alpha_j )_{j \geq 0}$ by
\[
    \alpha_j \equiv p^{-n} \BN(l_n + j\eulerphi(p^n))  \pmod{p^2} .
\]
Using Corollary \ref{corl-kummer-congr-irrpair}
with $r=3$, $\omega = \eulerphi(p^n)$, and $m=l_n+j \eulerphi(p^n) > 2$,
we have for $n(r-1) \geq 2$ that
\[
   \alpha_j - 3 \alpha_{j+1} + 3 \alpha_{j+2} - \alpha_{j+3} \equiv 0 \pmod{p^2} .
\]
Taking differences with $\beta_j = \alpha_{j+1} - \alpha_j$ yields
\[
   \beta_j - 2 \beta_{j+1} + \beta_{j+2} \equiv 0 \pmod{p^2} .
\]
By Lemma \ref{lemma-seq-a2} the sequence $( \beta_j )_{j \geq 0}$ is equidistant
$\mods{p^2}$. Proposition \ref{prop-irrpair-n} shows that
the sequence $( \alpha_j )_{j \geq 0}$ is equidistant $\mods{p}$.
By definition we have $\beta_j \equiv \Delta_n \modp{p}$.
Therefore, we can choose suitable $\gamma,\delta \in \ZZ$ so that
\[
   \alpha_{j+1} - \alpha_j \equiv \beta_j \equiv
     \Delta_n + p \, ( \gamma + j \delta ) \pmod{p^2} .
\]
This yields
\begin{eqnarray} \nonumber
   \alpha_j &\equiv& \alpha_0 + \sum_{\nu=0}^{j-1} (
     \Delta_n + p \, ( \gamma + \nu \delta ) ) \\
      &\equiv& \alpha_0 + j \Delta_n + j p \gamma +
        \binom{j}{2} p \delta \pmod{p^2} .  \label{eqn-loc-irr-delta-n-1}
\end{eqnarray}
From Proposition \ref{prop-irrpair-n} we have
\begin{equation} \label{eqn-loc-irr-delta-n-2}
  s \equiv -\alpha_0 \, \Delta_n^{-1} \pmod{p} , \quad 0 \leq s < p .
\end{equation}
With $l_{n+1} = l_n + s \eulerphi(p^n)$ we obtain the unique related irregular pair
$(p,l_{n+1}) \in \IRR_{n+1}$ of order $n+1$.
As a consequence of Lemma \ref{lemma-seq-a2}, we observe that
\[
   \alpha_{s+jp} \equiv 0 \pmod{p} .
\]
Thus we obtain a sequence $( \alpha'_j )_{j \geq 0}$ defined by
\[
   \alpha'_j \equiv \alpha_{s+jp} / p \equiv
     p^{-{(n+1)}} \BN(l_{n+1} + j \eulerphi(p^{n+1})) \pmod{p}
\]
which we can use to determine related irregular pairs of order $n+2$
using Proposition \ref{prop-irrpair-n}.
By definition we have
\[
    \Delta_{n+1} \equiv \alpha'_1 - \alpha'_0 \pmod{p} .
\]
It follows from \refeqn{eqn-loc-irr-delta-n-1} that
\begin{eqnarray*}
   p \, \Delta_{n+1} &\equiv& p \, (\alpha'_1 - \alpha'_0)
     \, \equiv \, \alpha_{s+p} - \alpha_s \\
       &\equiv& p \, \Delta_n + p^2 \, \gamma + p \, \delta \left(
         \binom{s+p}{2} - \binom{s}{2} \right) \\
       &\equiv& p \, \Delta_n \hspace*{33ex}  \pmod{p^2} ,
\end{eqnarray*}
since
\[
   \binom{s+p}{2} - \binom{s}{2} = \frac{1}{2} \, p \, (p+2s-1) .
\]
Finally, we obtain the proposed equation $\Delta_{n+1} = \Delta_n$.
\end{proof}

\begin{proofof}{Theorem~\ref{theor-delta-irr}}
Using Proposition \ref{prop-irr-delta-n} with induction on $n$ provides
\[
   \Delta_{(p,l_1)} = \Delta_{(p,l_2)} = \Delta_{(p,l_3)} = \ldots
\]
with exactly one related irregular pair of order $n$
\[
   (p,l_n) \in \IRR_n \quad \mbox{resp.} \quad (p,s_1,\ldots,s_n) \in \IRP_n .
\]
The latter pair is given by Definition \ref{def-pair-padic}.
Proposition \ref{prop-irrpair-n} shows for each step $n$ that
\begin{equation} \label{eqn-loc-delta-irr-1}
    l_{n+1} = l_n + s_n \, \eulerphi(p^n) , \quad  0 \leq s_n < p
\end{equation}
with a suitable integer $s_n$. This ensures that $l_1 \leq l_2 \leq l_3 \leq \ldots $ as
an increasing sequence $(l_j)_{j \geq 1}$. Clearly, this sequence is not eventually constant,
because $p^n \pdiv \BN(l_n)$ with $0 < |\BN(l_n)| < \infty$. Therefore
$\lim_{n \to \infty} l_n = \infty$. \smallskip

Starting with $(p,l_1) \in \IRR_1$ this provides a unique sequence $(l_j)_{j \geq 1}$.
If we have another irregular pair $(p,l'_1) \in \IRR_1$ with $\Delta_{(p,l'_1)} \neq 0$
and $(p,l_1) \neq (p,l'_1)$ then
\[
   l'_j \neq l_k \quad \mbox{for all\ } j,k \in \NN ,
\]
because $l'_1 \not\equiv l_1 \modp{\eulerphi(p)}$
and $l'_1 \equiv l'_j \modp{\eulerphi(p)}$ resp.\ $l_1 \equiv l_k \modp{\eulerphi(p)}$
by \refeqn{eqn-loc-delta-irr-1}. \smallskip

If $\Delta(p)=1$ then for each of the $i(p)$ irregular pairs $(p,l_{1,\nu})$,
$\nu=1,\ldots,i(p)$ there exists exactly one related irregular pair of higher order.
Finally, it follows that $i(p) = i_2(p) = i_3(p) = \ldots$ and so on.
\end{proofof}

\begin{proofof}{Theorem~\ref{theor-delta-irr-0}}
Clearly, the (non-) existence of related irregular pairs in the case (1) resp.\ (2)
is given by Proposition \ref{prop-irrpair-n} case (1) resp.\ (2).
Hence, we only have to show the remaining part of case (2).
In this case we have $p$ related irregular pairs
$(p,l_{n+1,\nu}) = (p,l_n+\nu \eulerphi(p^n)) \in \IRR_{n+1}$
of order $n+1$ with $\nu=0,\ldots,p-1$. Although $\Delta_{(p,l_n)}=0$ we can
use Proposition \ref{prop-irr-delta-n} by modifying essential steps.
We then have
\[
    \alpha_j \equiv \beta_j \equiv 0 \pmod{p} .
\]
Congruence \refeqn{eqn-loc-irr-delta-n-2} must be replaced by
\begin{equationname}{\ref{eqn-loc-irr-delta-n-2}'}
   s = 0,\ldots,p-1
\end{equationname}%
since $\Delta_{(p,l_n)}=0$ yields $p$ values of $s$. It follows that
\[
   \Delta_{(p,l_n)} = \Delta_{(p,l_{n+1,\nu})} = 0 \quad \mbox{for} \quad \nu=0,\ldots,p-1 .
     \vspace*{-4ex}
\]
\end{proofof}

\section{A $p$-adic view}
\setcounter{equation}{0}
\label{sect-padic}

Let $\ZZ_p$ be the ring of $p$-adic integers and $\QQ_p$ be the
field of $p$-adic numbers. The ultrametric absolute value $|\ |_p$
is defined by $|x|_p = p^{-\ord_p x}$ on $\QQ_p$.
Define $| \ |_\infty$ as the standard norm on $\QQ$.
For $n \geq 1$ we define $\psi_n: \ZZ_p \to \ZZ$ giving
the projection of $\ZZ_p$ onto the set $[0,p^n) \cap \ZZ$ so that
for $x \in \ZZ_p$ we have $x-\psi_n(x) \in p^n \ZZ_p$ where
$0 \leq \psi_n(x) < p^n$.
We denote $\PP$ as the set of the rational primes.
Now, we shall use \refeqn{eqn-zeta-1-n} to reformulate our results.
Let $(p,l) \in \IRR_1$ then
\[
    \Delta_{(p,\,l)} \, \equiv \,
       \frac{\zeta(1-l) \, - \, \zeta(1-(l+p-1))}{p} \pmod{p}
\]
with $0 \leq \Delta_{(p,\,l)} < p$.

\begin{theorem} \label{theor-zeta-padic}
Let $(p,l) \in \IRR_1$ with $\Delta_{(p,\,l)} \neq 0$.
We define the sequence $(l_n)_{n \geq 1}$ recursively
by $l_1 = l$ and, for $n \geq 1$, by
\[
   l_{n+1} \,=\, l_n + \eulerphi(p) \, \psi_n \!
       \left( \frac{\zeta(1-l_n)}{p \, \Delta_{(p,\,l)}} \right)
   \,=\, l_n + \eulerphi(p^n) \, \psi_1 \!
     \left( \frac{\zeta(1-l_n)}{p^n \, \Delta_{(p,\,l)}} \right) \!.
\]
Then we have $\zeta(1-l_n) \in p^n \ZZ_p$ and consequently
\[
   \lim_{n \to \infty} |\zeta(1-l_n)|_p = 0
    \quad \mbox{with} \quad l_n \to \infty .
\]
\end{theorem}

\begin{proof}
We rewrite our results using $\zeta(1-l_n) = - \BN(l_n)$.
Theorem \ref{theor-delta-irr} provides one, and only one,
sequence $(l_n)_{n \geq 1}$
with $l_1 = l$ and $(p,l_n) \in \IRR_n$ with $l_n \to \infty$.
This implies the $p$-adic convergence $\zeta(1-l_n) \to 0$.
Additionally, from Proposition \ref{prop-irrpair-n} we have
\[
   l_{n+1} = l_n + s \, \eulerphi(p^n)
\]
for each step where
\[
   s \equiv - p^{-n} \BN(l_n) \, \Delta_{(p,\,l)}^{-1} \pmod{p}
     , \quad 0 \leq s < p .
\]
Rewriting the last congruence yields
\[
   s = \psi_1 \! \left( \frac{\zeta(1-l_n)}{p^n \, \Delta_{(p,\,l)}} \right) \!.
\]
The rest follows by $\psi_n( a \, p^{n-1} ) =  p^{n-1} \, \psi_1(a)$ for $a \in \ZZ_p$.
\end{proof}

These results can be also applied to the so-called $p$-adic zeta functions
which were originally defined by T.~Kubota and H.~W.~Leopoldt \cite{leopoldt64} in 1964;
for a detailed theory see Koblitz \cite[Chapter II]{koblitz96}.
We introduce some definitions.

\begin{defin} \label{def-zeta-padic}
Let $p$ be a prime. For $n \geq 1$ define
\[
   \zeta_p(1-n) \, := \, (1-p^{n-1}) \, \zeta(1-n)
       \, = \, -(1-p^{n-1}) \, \BN(n) .
\]
Define the \textsl{$p$-adic zeta function}
for $p \geq 5$ and a fixed $s_1 \in \{ 2,4,\ldots,p-3 \}$ by
\[
   \zeta_{p,\,s_1}:\,\, \ZZ_p \to \ZZ_p , \qquad
     \zeta_{p,\,s_1} (s) := \lim_{t_\nu \to s} \,
       \zeta_p \big( 1-(s_1 +(p-1) t_\nu) \big)
\]
resp.\ for $p \geq 2$ and $s_1 = 0$ by
\[
   \zeta_{p,\,0}:\,\, \ZZ_p \setminus \{0\} \to \QQ_p , \qquad
     \zeta_{p,\,0} (s) := \lim_{t_\nu \to s} \,
       \zeta_p \big( 1-(p-1) t_\nu \big)
\]
for $p$-adic integers $s$ by taking any sequence $(t_\nu)_{\nu \geq 1}$ of nonnegative
integers resp.\ positive integers in the case $s_1 = 0$ which $p$-adically converges to $s$.
\end{defin}

\begin{remark} \label{rem-padic-zeta}
Let $p \geq 2$ and $s_1 \geq 0$.
The $p$-adic zeta function $\zeta_{p,\,s_1} (s)$
interpolates the zeta function $\zeta_p(1-n)$ at nonnegative integer values $s$ by
\[
   \zeta_{p,\,s_1}(s) = \zeta_p(1-n)
\]
where $n \equiv s_1 \modp{p-1}$ and $n=s_1+(p-1)s$.

Let $p \geq 5$ and $s_1 \neq 0$.
The Kummer congruences \refeqn{eqn-kummer-congr} state for $r \geq 1$ that
\[
   \zeta_{p,\,s_1}(s) \equiv \zeta_{p,\,s_1}(s') \pmod{p^r}
\]
when $s \equiv s' \modp{p^{r-1}}$ for nonnegative integers $s$ and $s'$.
Since $\ZZ$ is dense in $\ZZ_p$,
the $p$-adic zeta function $\zeta_{p,\,s_1}$ extends uniquely,
by means of the interpolation property and the Kummer congruences,
to a continuous function on $\ZZ_p$; see \cite[Thm.~8, p.~46]{koblitz96}.
\end{remark}

\begin{defin} \label{def-char-padic}
Let $(p,l) \in \IRR_1$ with $\Delta_{(p,\,l)} \neq 0$.
Define a \textsl{characteristic} $p$-adic integer
which contains all information about irregular pairs of higher order by
\[
   \chi_{(p,\,l)} = \sum_{\nu \geq 0} s_{\nu+2} \, p^\nu \,\, \in \ZZ_p
\]
where $(s_\nu)_{\nu \geq 1}$ is the sequence defined by Theorem \ref{theor-delta-irr}
and $l=s_1$.
\end{defin}

\begin{lemma} \label{lemma-ordp_bn_chi}
Let $(p,l) \in \IRR_1$ with $\Delta_{(p,\,l)} \neq 0$.
Let $r = \ord_p \BN(l)$ and $(p,s_1,\ldots,s_{r+1}) \in \IRP_{r+1}$ be
the related irregular pair of order $r+1$. Then
\[
   s_{r+1} \, \Delta_{(p,\,l)} \equiv - p^{-r} \, \BN(l) \pmod{p}
\]
with $s_1=l$, $s_\nu=0$ for $\nu=2,\ldots,r$, $s_{r+1} \neq 0$.
If $r=1$ then $\chi_{(p,\,l)} \in \ZZ_p^*$, otherwise
$\chi_{(p,\,l)} \in p^{r-1} \ZZ_p$ for $r \geq 2$.
\end{lemma}

\begin{proof}
Since $r = \ord_p \,\BN(l) \geq 1$, we have $(p,l) \in \IRR_\nu$ for all $\nu=1,\ldots,r$.
Then Proposition \ref{prop-irrpair-n} and Theorem \ref{theor-delta-irr} provide
\[
   s \equiv - p^{-r} \BN(l) \, \Delta_{(p,\,l)}^{-1} \pmod{p} ,
     \qquad (p,l+s \, \eulerphi(p^r)) \in \IRR_{r+1}
\]
with $0 \leq s < p$. We have $s \neq 0$ since $\ord_p\,(p^{-r} \BN(l)) = 0$.
Let $(p,s_1,\ldots,s_{r+1}) \in \IRP_{r+1}$ be the related irregular pair of order $r+1$.
Then we see that $s_1 = l$, $s_{r+1} = s \neq 0$, and $s_\nu=0$ for $\nu=2,\ldots,r$.
By Definition \ref{def-char-padic} we have
\[
   \chi_{(p,\,l)} = s_2 + s_3 \, p + \ldots + s_{r+1} \, p^{r-1} + \ldots \,.
\]
Case $r=1$ yields $s_2 \neq 0$ and $\chi_{(p,\,l)} \in \ZZ_p^*$, otherwise
case $r \geq 2$ implies that $\ord_p \chi_{(p,\,l)} = r-1$ and
$\chi_{(p,\,l)} \in p^{r-1} \ZZ_p$.
\end{proof}

\begin{theorem} \label{theor-zeta-padic-zero}
Let $(p,l) \in \IRR_1$ with $\Delta_{(p,\,l)} \neq 0$.
The $p$-adic zeta function $\zeta_{p,\,l}(s)$ has a unique zero
at $s=\chi_{(p,\,l)}$.
\end{theorem}

\begin{proof}
From Theorem \ref{theor-delta-irr} we have a sequence
$(l_n)_{n \geq 1}$ with $l=l_1$.
In view of Definition \ref{def-zeta-padic},
Theorem \ref{theor-zeta-padic} also states that
\[
   \lim_{n \to \infty} |\zeta_p(1-l_n)|_p = 0
      \quad \mbox{with} \quad l_n \to \infty .
\]
We can transfer this result to the $p$-adic zeta function $\zeta_{p,\,l}$
by the interpolation property. We see that $p$-adically
\[
    \lim_{n \to \infty} l_n = l+(p-1)\chi_{(p,\,l)} .
\]
Since the function $\zeta_{p,\,l}$ is continuous, the
$p$-adic integer $\chi_{(p,\,l)}$ is a zero of $\zeta_{p,\,l}$.
We shall show that this zero is unique.
Assume that $\zeta_{p,\,l}(\xi) = 0$ with some $\xi \in \ZZ_p$.
We can use the arguments given above in the opposite direction.
Since $\zeta_{p,\,l}$ is continuous, there exists a sequence $(l'_n)_{n \geq 1}$
of positive integers with
\[
    \lim_{n \to \infty} l'_n = l+(p-1) \xi \quad \mbox{and} \quad
      \lim_{n \to \infty} |\zeta_p(1-l_n)|_p = 0 .
\]
We can choose a subsequence $(l''_n)_{n \geq 1}$ of $(l'_n)_{n \geq 1}$
such that $\zeta_p(1-l''_n) \in p^n \ZZ_p$.
By use of the Kummer congruences, we construct the
sequence $(\tilde{l}_n)_{n \geq 1}$ where $\tilde{l}_n \equiv l''_n \modp{\eulerphi(p^n)}$
with $l \leq \tilde{l}_n < \eulerphi(p^n)$. Now we have $(p,\tilde{l}_n) \in \IRR_n$ for all $n \geq 1$.
Since $l=l_1=\tilde{l}_1$, Theorem \ref{theor-delta-irr} shows that
$(l_n)_{n \geq 1}$ = $(\tilde{l}_n)_{n \geq 1}$. This implies that $\xi = \chi_{(p,\,l)}$.
\end{proof}

\begin{remark}
Let $(p,l) \in \IRR_1$ with $\Delta_{(p,\,l)} \neq 0$.
The related irregular pair $(p,s_1,\ldots,s_r)$ $\in$ $\IRP_r$
is a $p$-adic approximation
of the zero $\chi_{(p,\,l)}$ of the $p$-adic zeta function $\zeta_{p,\,l}$.
For the first irregular primes 37, 59, and 67 elements of $\IRP_{100}$ were
calculated in \cite[pp.~127--128]{kellner02}. These results are
reprinted in Table \ref{tab-irrpair-ord-100}.
By Lemma \ref{lemma-ordp_bn_chi} the statement $p^2 \notdiv B_l$
is equivalent to the fact that $\chi_{(p,\,l)}$ is a unit in $\ZZ_p$.
\end{remark}
\smallskip

From now on, we shall assume the $\Delta$-Conjecture.
We shall see that the zeros $\chi_{(p,\,l)}$ play an important role
in the representation of the Riemann zeta function at odd negative integer arguments.

\begin{theorem} \label{theor-zeta-prod-psi}
Let $n$ be an even positive integer.
Under the assumption of the $\Delta$-Conjecture we have
\[
   \zeta(1-n) = (-1)^{\frac{n}{2}} \, \prod_{p-1 \pdiv n} \frac{|n|_p}{p} \,
     \prod_{(p,l)\in\IRR_1 \atop l \equiv n \modp{p-1}}
       \frac{p}{\,|\chi_{(p,\,l)} - \frac{n-l}{p-1}|_p} .
\]
\end{theorem}

\begin{proof}
Since both products above have only positive terms, the sign follows by \refeqn{eqn-zeta-1-n}.
The first product describes the denominator of $\zeta(1-n)$ which is
a consequence of \refeqn{eqn-clausen-von-staudt} and \refeqn{eqn-adams-triv-div}.
We have to show that the second product describes the unsigned numerator of $\zeta(1-n)$
which only consists of powers of irregular primes.
Let $p$ be an irregular prime divisor of $\zeta(1-n)$.
From Remark \ref{rem-ordn-irrpair} we have
\[
   \ord_p \BN(n) = r \quad \Longrightarrow \quad n \equiv l_r \pmod{\eulerphi(p^r)}
     \quad \mbox{with} \quad (p,l_r) \in \IRR_r .
\]
The irregular pair $(p,l_r)$ of order $r$ is related to some irregular pair
$(p,l) \in \IRR_1$ with $\Delta_{(p,\,l)} \neq 0$
where $l \equiv l_r \equiv n \modp{\eulerphi(p)}$.
We also have by Definition \ref{def-char-padic} that
\[
   n \equiv l_r \equiv l+(p-1) \, \chi_{(p,\,l)} \pmod{(p-1)p^{r-1}\ZZ_p}
\]
and equally by reduction that
\[
   \frac{n-l}{p-1} \equiv \chi_{(p,\,l)} \pmod{p^{r-1}\ZZ_p} .
\]
The last congruence is not valid $\mods{p^r\ZZ_p}$ by construction. Therefore
\[
    \left| \chi_{(p,\,l)} - \frac{n-l}{p-1} \right|_p = p^{-(r-1)}
\]
which provides, with an additional factor $p$, the desired $p$-power
in the second product. Considering all irregular primes $p$ which can appear,
the second product equals the numerator of $\zeta(1-n)$ without sign.
\end{proof}

With some technical definitions we can combine both products of the theorem above.
This yields a more accessible representation of the Riemann zeta function by means of
$p$-adic analysis.

\begin{theorem} \label{theor-zeta-prod-psi-0}
Define $\Psi_0 = \IRR_1 \cup ( \PP \times \{0\} )$
and set $\chi_{(p,\,0)}=0$ for all $p \in \PP$.
Define $\rho(l)=1-2\sign(l)=\pm1$ for $l \geq 0$.
Let $n$ be an even positive integer.
Under the assumption of the $\Delta$-Conjecture we have
\[
   \zeta(1-n) = (-1)^{\frac{n}{2}} \!\!\!
     \prod_{(p,l)\in\Psi_0 \atop l \equiv n \modp{p-1}}
       \left( \frac{\,|\chi_{(p,\,l)} - \frac{n-l}{p-1}|_p}{p} \right)^{\rho(l)} .
\]
\end{theorem}

\begin{proof}
We only have to consider case $l=0$. Then we have $p-1 \pdiv n$, $\rho(0)=1$,
and $|\chi_{(p,\,0)} - \frac{n}{p-1}|_p = |n|_p$. The other case $l > 0$ is
already covered by Theorem \ref{theor-zeta-prod-psi}.
\end{proof}
\smallskip

We shall give an interpretation of this formula above in Remark \ref{rem-zeta-prod}
by means of $p$-adic zeta functions. This generalization
shows the significance to prove the $\Delta$-Conjecture at all.

\begin{theorem} \label{theor-zeta-padic-delta-kummer}
Let $(p,l) \in \IRR_1$ with $\Delta_{(p,\,l)} \neq 0$.
Let $s, t \in \ZZ_p$. Then a strong version of the Kummer congruences holds that
\[
   |\zeta_{p,\,l} (s) - \zeta_{p,\,l} (t)|_p = |p \, ( s - t ) |_p .
\]
Moreover
\[
   \frac{\zeta_{p,\,l} (s) - \zeta_{p,\,l} (t)}{p \, ( s - t )}
      \equiv - \Delta_{(p,\,l)} \pmod{p\ZZ_p}
         \quad \mbox{for} \quad s \neq t
\]
and
\[
   \zeta'_{p,\,l} (s) \equiv - p \, \Delta_{(p,\,l)} \pmod{p^2\ZZ_p} .
\]
\end{theorem}

Thus, $\Delta_{(p,\,l)}$ is closely associated with
the $p$-adic zeta function $\zeta_{p,\,l}$ in the
nonsingular case $\Delta_{(p,\,l)} \neq 0$. We will prove this theorem later.

\begin{corl} \label{corl-zeta-padic-zero}
Let $(p,l) \in \IRR_1$ with $\Delta_{(p,\,l)} \neq 0$.
The $p$-adic zeta function $\zeta_{p,\,l}(s)$ has a simple zero
at $s=\chi_{(p,\,l)}$. Moreover, for $s \in \ZZ_p$,
\[
   \zeta_{p,\,l}(s) = p \, ( s - \chi_{(p,\,l)} ) \, \zeta^\star_{p,\,l}(s)
\]
where $\zeta^\star_{p,\,l} (s)$ is a continuous function on $\ZZ_p$ with
$\zeta^\star_{p,\,l} (s) \equiv - \Delta_{(p,\,l)} \pmod{p \ZZ_p}$.
Consequently
\[
   |\zeta_{p,\,l}(s)|_p = |p \, ( s - \chi_{(p,\,l)} ) |_p .
\]
\end{corl}

\begin{proof}
Theorem \ref{theor-zeta-padic-zero} shows that  $\zeta_{p,\,l}(s)$ has
a unique zero at $s = \chi_{(p,\,l)}$.
We can use Theorem \ref{theor-zeta-padic-delta-kummer} to define
\[
   \zeta^\star_{p,\,l}(s) = \left\{
     \begin{array}{cl} \displaystyle
       \frac{\zeta_{p,\,l}(s)}{p(s-\chi_{(p,\,l)})} , & s \neq \chi_{(p,\,l)} , \\\\
       \displaystyle
       \frac{\zeta'_{p,\,l}(s)}{p} , & s = \chi_{(p,\,l)} . \\
     \end{array} \right.
\]
Theorem \ref{theor-zeta-padic-delta-kummer} implies
that $\zeta^\star_{p,\,l} (s) \equiv - \Delta_{(p,\,l)} \pmod{p \ZZ_p}$
for all $s \in \ZZ_p$. Hence $\zeta^\star_{p,\,l} (s)$ has
no zeros and consequently $\zeta_{p,\,l}(s)$ has a simple zero at
$s = \chi_{(p,\,l)}$.
Since $\zeta_{p,\,l}(s)$ is continuous on $\ZZ_p$ and
$\zeta'_{p,\,l}(s)$ exists at $s=\chi_{(p,\,l)}$,
the function $\zeta^\star_{p,\,l}(s)$ is also continuous on $\ZZ_p$.
Finally we obtain
$\zeta_{p,\,l}(s) = p \, ( s - \chi_{(p,\,l)} ) \, \zeta^\star_{p,\,l}(s)$
and $|\zeta_{p,\,l}(s)|_p = |p \, ( s - \chi_{(p,\,l)} ) |_p \, |\zeta^\star_{p,\,l}(s)|_p
= |p \, ( s - \chi_{(p,\,l)} ) |_p$.
\end{proof}

\begin{defin} \label{def-zeta-padic-0}
Let $p$ be a prime. Define
\[
   \zeta^\star_{p,\,0} : \ZZ_p \to \ZZ_p ,
     \qquad \zeta^\star_{p,\,0}(s) :=
       - \lim_{t_\nu \to s} (1-p^{t_\nu(p-1)-1}) \frac{p B_{t_\nu(p-1)}}{p-1}
\]
for $p$-adic integers $s$ by taking any sequence $(t_\nu)_{\nu \geq 1}$ of nonnegative
integers which $p$-adically converges to $s$.
\end{defin}

\begin{prop} \label{prop-zeta-padic-p-0}
Let $p$ be a prime. The function $\zeta^\star_{p,\,0}(s)$ is continuous on $\ZZ_p \setminus \{0\}$
and satisfies the following properties:
\[
   \zeta^\star_{p,\,0}(0) = -1 , \quad
     \zeta^\star_{p,\,0}(s) = ps \, \zeta_{p,\,0}(s) ,
     \quad s \in \ZZ_p \setminus \{0\}
\]
and
\[
   \lim_{s \to 0} |ps \, \zeta_{p,\,0}(s) - \zeta^\star_{p,\,0}(0)|_p < 1 .
\]
Moreover $\zeta^\star_{p,\,0}(s) \equiv -1 \pmod{p \ZZ_p}$
for $s \in \ZZ_p$ in case $p > 2$ and $s \in p\ZZ_p$ in case $p = 2$.
Additionally, if $p=2$ then $\zeta^\star_{p,\,0}(s) = 0$ for $s \equiv 1 \pmod{p \ZZ_p}$.
\end{prop}

\begin{proof}
From Definition \ref{def-zeta-padic} and Definition \ref{def-zeta-padic-0} we have
for $s \in \ZZ_p \setminus \{0\}$ and any sequence $(t_\nu)_{\nu \geq 1}$ of positive
integers which $p$-adically converges to $s$ that
\begin{eqnarray*}
   ps \, \zeta_{p,\,0}(s)
     &=& - \lim_{t_\nu \to s} (1-p^{t_\nu(p-1)-1}) ps \frac{B_{t_\nu(p-1)}}{t_\nu(p-1)} \\
     &=& - \lim_{t_\nu \to s} (1-p^{t_\nu(p-1)-1}) \frac{p B_{t_\nu(p-1)}}{p-1}
       \quad = \quad \zeta^\star_{p,\,0}(s) .
\end{eqnarray*}
Moreover we have
\[
   \zeta^\star_{p,\,0}(0) = -(1-p^{-1}) \, p B_0 / (p-1) = -1 .
\]
By Clausen-von Staudt \refeqn{eqn-clausen-von-staudt}
we obtain $\zeta^\star_{p,\,0}(s) \equiv -1 \pmod{p \ZZ_p}$
for $s \in \ZZ_p$ in case $p > 2$ and $s \in p\ZZ_p$ in case $p = 2$.
Additionally we have $\zeta^\star_{2,\,0}(s) = 0$ for $s \in 1 + 2\ZZ_2$,
since $\zeta^\star_{2,\,0}(1) = -(1-2^0) \, 2 B_1  = 0$ and
$B_n=0$ for all odd integers $n > 1$.
We use the fact that $\zeta_{p,\,0}: \ZZ_p \to \QQ_p$ is a continuous function on
$\ZZ_p \setminus \{0\}$; see \cite[Thm.~8, p.~46]{koblitz96}.
Hence $\zeta^\star_{p,\,0}(s) = ps \, \zeta_{p,\,0}(s)$
is also continuous on $\ZZ_p \setminus \{0\}$.
It remains to show that
\[
   \lim_{s \to 0} |ps \, \zeta_{p,\,0}(s) - \zeta^\star_{p,\,0}(0)|_p < 1 .
\]
We can choose a zero sequence $(s_\nu)_{\nu \geq 1}$
where its elements are arbitrary close to 0, say $0 < |s_\nu|_p < p^{-r}$
with some positive integer $r$. Then we deduce that
\[
   ps_\nu \, \zeta_{p,\,0}(s_\nu) - \zeta^\star_{p,\,0}(0)
     \equiv \zeta^\star_{p,\,0}(s_\nu) +1 \equiv 0 \pmod{p \ZZ_p}
\]
for all $\nu \geq 1$. This implies, together with the continuity of
$ps \, \zeta_{p,\,0}(s)$, the estimate above.
\end{proof}

\begin{corl} \label{corl-zeta-padic-pole}
Let $p$ be a prime. The $p$-adic zeta function $\zeta_{p,\,0}(s)$
has a simple pole at $s=0$. Moreover, for $s \in \ZZ_p \setminus \{0\}$,
\[
   \zeta_{p,\,0}(s) = \zeta^\star_{p,\,0}(s) / ( ps ) .
\]
Consequently
\[
   |\zeta_{p,\,0}(s)|_p = |ps|_p^{-1}
\]
for $s \in \ZZ_p \setminus \{0\}$ in case $p > 2$ and
$s \in p\ZZ_p \setminus \{0\}$ in case $p = 2$.
\end{corl}

\begin{proof}
By Proposition \ref{prop-zeta-padic-p-0} we can write
$\zeta_{p,\,0}(s) = \zeta^\star_{p,\,0}(s) / ( ps )$
for $s \in \ZZ_p \setminus \{0\}$.
Moreover $\lim_{s \to 0} ps \, \zeta_{p,\,0}(s) = \xi$ with
$|\xi-\zeta^\star_{p,\,0}(0)|_p < 1$.
This implies that the limit $\lim_{s \to 0} \zeta_{p,\,0}(s)$ does not exist
and that $\zeta_{p,\,0}(s)$ has a simple pole at $s=0$.
We have $|\zeta^\star_{p,\,0}(s)|_p = 1$ for $s \in \ZZ_p \setminus \{0\}$ in case $p > 2$ and
$s \in p\ZZ_p \setminus \{0\}$ in case $p = 2$.
In these cases we then obtain $|\zeta_{p,\,0}(s)|_p = |ps|_p^{-1}$.
\end{proof}

\begin{remark}
One can even show further that $\zeta^\star_{p,\,0}(s)$ is continuous on $\ZZ_p$ where
$\lim_{s \to 0} \zeta^\star_{p,\,0}(s) = \zeta^\star_{p,\,0}(0)$; moreover
$\zeta^\star_{p,\,0}(s)$ satisfies the Kummer congruences.
This can be derived, e.g., by means of $p$-adic integration, see Koblitz
\cite[pp.~42--46]{koblitz96}, or by using certain congruences
of Carlitz \cite{carlitz60}.
\end{remark}

\begin{theorem} \label{theor-zeta-pole-zero}
Define $\Psi_0 = \IRR_1 \cup ( \PP \times \{0\} )$
and set $\chi_{(p,\,0)}=0$ for all $p \in \PP$.
Define $s_{p,\,l}(n) = (n-l)/(p-1)$.
Let $n$ be an even positive integer. Then
\[
   |\zeta(1-n)|_\infty = \prod_{p \in \PP \atop l \equiv n \modp{p-1}}
     \hspace*{-3ex} |\zeta_{p,\,l}(s_{p,\,l}(n))|_p^{-1}
   = \prod_{(p,l)\in \Psi_0 \atop l \equiv n \modp{p-1}}
     \hspace*{-3ex} |\zeta_{p,\,l}(s_{p,\,l}(n))|_p^{-1} .
\]
Under the assumption of the $\Delta$-Conjecture we have
\[
   |\zeta(1-n)|_\infty
   = \prod\limits_{(p,l)\in \IRR_1 \atop l \equiv n \modp{p-1}}
     \hspace*{-3ex} |p(s_{p,\,l}(n)-\chi_{(p,\,l)})|_p^{-1}
     \Big/ \prod_{p-1 \pdiv n}
     |p(s_{p,\,0}(n)-\chi_{(p,\,0)})|_p^{-1} .
\]
\end{theorem}

\begin{proof}
Since $n$ is an even positive integer, the product formula states that
\[
   \prod_{p \, \in \, \PP \, \cup \, \{\infty\}} |\zeta(1-n)|_p = 1 .
\]
By Definition \ref{def-zeta-padic} we have
$|\zeta(1-n)|_p = |\zeta_{p,\,l}(s)|_p$
where $l \equiv n \modp{p-1}$ with $0 \leq l < p-1$ and
$s=s_{p,\,l}(n)=(n-l)/(p-1)$. Thus
\[
   |\zeta(1-n)|_\infty^{-1} = \prod_{p \, \in \, \PP} |\zeta(1-n)|_p
   = \prod_{p \in \PP \atop l \equiv n \modp{p-1}}
     \hspace*{-3ex} |\zeta_{p,\,l}(s_{p,\,l}(n))|_p \, .
\]
We have $|\zeta_{p,\,l}(s)|_p=1$ for $s \in \ZZ_p$
if $(p,l) \notin \IRR_1$ and $l \neq 0$.
Hence
\[
   \prod_{p \in \PP \atop l \equiv n \modp{p-1}}
     \hspace*{-3ex} |\zeta_{p,\,l}(s_{p,\,l}(n))|_p
   = \prod_{(p,l)\in \Psi_0 \atop l \equiv n \modp{p-1}}
     \hspace*{-3ex} |\zeta_{p,\,l}(s_{p,\,l}(n))|_p \, .
\]
Next we assume the $\Delta$-Conjecture. By Corollary \ref{corl-zeta-padic-zero} we
then have for $(p,l) \in \IRR_1$ that
$|\zeta_{p,\,l}(s)|_p = |p \, ( s - \chi_{(p,\,l)} ) |_p$ for $s \in \ZZ_p$.
Without any assumption, Corollary \ref{corl-zeta-padic-pole} shows that
$|\zeta_{p,\,0}(s)|_p = |ps|_p^{-1} = |p(s-\chi_{(p,\,0)})|_p^{-1}$
for $s \in \ZZ_p \setminus \{0\}$ in case $p > 2$ and
$s \in p\ZZ_p \setminus \{0\}$ in case $p = 2$.
Since $n$ is even and $s_{2,\,0}(n)=n$, we finally obtain
\[
   \prod_{(p,l)\in \Psi_0 \atop l \equiv n \modp{p-1}}
     \hspace*{-3ex} |\zeta_{p,\,l}(s_{p,\,l}(n))|_p
   = \hspace*{-3.5ex} \prod\limits_{(p,l)\in \IRR_1 \atop l \equiv n \modp{p-1}}
     \hspace*{-3ex} |p(s_{p,\,l}(n)-\chi_{(p,\,l)})|_p
     \prod_{p-1 \pdiv n} |p(s_{p,\,0}(n)-\chi_{(p,\,0)})|_p^{-1} . \vspace*{-2ex}
\]
\end{proof}

\begin{remark} \label{rem-zeta-prod}
Assuming the $\Delta$-Conjecture,
the numerator of $|\zeta(1-n)|_\infty$ with even $n > 0$ is essentially described
by simple zeros $\chi_{(p,\,l)}$ of $p$-adic zeta functions $\zeta_{p,\,l}(s)$
where $l \equiv n \modp{p-1}$
and $(p,l) \in \IRR_1$. More precisely, by the variable substitution $s=(n-l)/(p-1)$,
the term
\[
   |p(s-\chi_{(p,\,l)})|_p^{-1}
\]
is equal to $p^r$ for some suitable $r > 0$ which is the $p$ power factor
in the numerator of $|\zeta(1-n)|_\infty$. The term above is induced by
\[
   |\zeta_{p,\,l}(s)|_p = |p(s-\chi_{(p,\,l)})|_p = p^{-1} |(s-\chi_{(p,\,l)})|_p \, .
\]
Since $\chi_{(p,\,l)}$ is a simple zero of $\zeta_{p,\,l}(s)$,
$|\zeta_{p,\,l}(s)|_p$ is mainly determined by a linear factor
which one can also interpret as a distance between $s$ and $\chi_{(p,\,l)}$.

Similar arguments can be applied to the denominator of $|\zeta(1-n)|_\infty$.
Without any assumption, the denominator of $|\zeta(1-n)|_\infty$ is essentially described
by simple poles $\chi_{(p,\,0)}$ of $p$-adic zeta functions
$\zeta_{p,\,0}(s)$ where $p-1 \pdiv n$. Again, the term
\[
   |p(s-\chi_{(p,\,0)})|_p^{-1}
\]
is equal to $p^r$ for some suitable $r > 0$ which is the $p$ power factor
in the denominator of $|\zeta(1-n)|_\infty$. This term is induced by
\[
   |\zeta_{p,\,0}(s)|_p = |p(s-\chi_{(p,\,0)})|_p^{-1} = p \, |(s-\chi_{(p,\,0)})|_p^{-1}
\]
where $\zeta_{p,\,0}(s)$ has a simple pole at $s = \chi_{(p,\,0)}$.
\end{remark}
\smallskip

Now we shall make some preparations to give later a proof of
Theorem \ref{theor-zeta-padic-delta-kummer}.

\begin{defin} \label{def-diff-oper}
Define the linear finite-difference operator $\Dop$ and its powers by
\[
   \Dop^r f(s) = \sum_{\nu=0}^r \binom{r}{\nu} (-1)^{\nu} f(s+\nu)
\]
for $r \geq 0$ and any function $f: \, \ZZ_p \to \ZZ_p$.
The series
\[
   f(s) = \sum_{\nu \geq 0} a_\nu \binom{s}{\nu}
\]
with coefficients $a_\nu \in \ZZ_p$ where $|a_\nu|_p \to 0$ is called a \textsl{Mahler series}
which defines a continuous function $f: \, \ZZ_p \to \ZZ_p$.
\end{defin}

The following theorem of Mahler shows
that the converse also holds, see \cite[Thm.~1, p.~173]{robert00}.
Note that the sign $(-1)^\nu$ depends on the definition of $\Dop$.

\begin{theorem}[Mahler] \label{theor-mahler}
Let $f: \, \ZZ_p \to \ZZ_p$ be a continuous function.
Then $f$ has a Mahler series
\[
   f(s) = \sum_{\nu \geq 0} a_\nu \binom{s}{\nu}
\]
where the coefficients $a_\nu$ are given by
$a_\nu = (-1)^\nu \, \Dop^\nu f(0)$ and $|a_\nu|_p \to 0$ holds.
\end{theorem}

\begin{defin} \label{def-zeta-padic-ordn}
Let $(p,l) \in \IRR_1$ with $\Delta_{(p,\,l)} \neq 0$.
For $n \geq 1$ define the $p$-adic zeta function of order $n$ by
\[
   \zeta_{p,\,l,\,n} (s) = p^{-n} \, \zeta_{p,\,l}
     \left( \psi_{n-1}( \chi_{(p,\,l)} ) + p^{n-1} s \right) \!,
       \quad s \in \ZZ_p
\]
which for $s=0$ corresponds to the related irregular pair of order $n$.
\end{defin}

\begin{prop} \label{prop-zeta-padic-ordn-diff}
Let $(p,l) \in \IRR_1$ with $\Delta_{(p,\,l)} \neq 0$.
For positive integers $n$, $r$ we have
\begin{equation} \label{eqn-loc-zeta-padic-ordn-diff-1}
   \Dop^{r+1} \zeta_{p,\,l,\,n} (s) \equiv 0 \pmod{p^{nr} \ZZ_p} ,
     \quad s \in \ZZ_p
\end{equation}
and
\begin{equation} \label{eqn-loc-zeta-padic-ordn-diff-2}
   \zeta_{p,\,l,\,n} (s) \equiv \zeta_{p,\,l,\,n} (t) \pmod{p^r \ZZ_p} ,
     \quad s, t \in \ZZ_p
\end{equation}
when $s \equiv t \modp{p^r \ZZ_p}$.
\end{prop}

\begin{proof}
First assume that $s \in \NN_0$.
In analogy to Corollary \ref{corl-kummer-congr-irrpair},
we have to modify Theorem \ref{theor-kummer-congr-general} in a similar way,
since $\zeta_{p,\,l,\,n} (0)$ corresponds to the related irregular pair of order $n$.
Rewriting \refeqn{eqn-kummer-congr-general} in the case $k=1$ gives
\refeqn{eqn-loc-zeta-padic-ordn-diff-1}.
Since $\ZZ$ is dense in $\ZZ_p$, we can extend
\refeqn{eqn-loc-zeta-padic-ordn-diff-1} to values $s \in \ZZ_p$
by means of the interpolation property of $\zeta_{p,\,l}$ resp.\ $\zeta_{p,\,l,\,n}$.
By the same arguments the Kummer congruences,
given in Remark \ref{rem-padic-zeta}, can be extended to values in $\ZZ_p$.
Let $s \equiv t \modp{p^r \ZZ_p}$ and
write $s' = \psi_{n-1}( \chi_{(p,\,l)} ) + p^{n-1} s$ and
$t' = \psi_{n-1}( \chi_{(p,\,l)} ) + p^{n-1} t$.
Then $s' \equiv t' \modp{p^{r+n-1} \ZZ_p}$ and therefore
$\zeta_{p,\,l} (s') \equiv \zeta_{p,\,l} (t') \modp{p^{r+n} \ZZ_p}$.
By Definition \ref{def-zeta-padic-ordn}
this gives \refeqn{eqn-loc-zeta-padic-ordn-diff-2}.
\end{proof}

\begin{prop} \label{prop-zeta-padic-ordn-exp}
Let $(p,l) \in \IRR_1$ with $\Delta_{(p,\,l)} \neq 0$.
For $n \geq 1$ we have
\[
   \zeta_{p,\,l,\,n} (s) \equiv \Delta_{(p,\,l)} \, ( s_{n+1} - s ) \pmod{p \ZZ_p} ,
       \quad s \in \ZZ_p
\]
where $s_{n+1}$ is defined by $\chi_{(p,\,l)} = s_2 + s_3 \, p + \ldots$ \,.
There exists the Mahler expansion
\[
   \zeta_{p,\,l,\,n}(s) = \zeta_{p,\,l,\,n}(0) +
     \sum_{\nu \geq 1} p^{n(\nu-1)} z_\nu \binom{s}{\nu} ,
        \quad s \in \ZZ_p
\]
with $z_\nu \in \ZZ_p$
where $z_\nu = (-1)^\nu p^{-n(\nu-1)} \, \Dop^\nu \zeta_{p,\,l,\,n}(0)$
and $z_1 \equiv -\Delta_{(p,\,l)} \modp{p \ZZ_p}$.
\end{prop}

\begin{proof}
Proposition \ref{prop-irrpair-n} also works with
$\alpha_j \equiv \zeta_{p,\,l,\,n}(j) \modp{p \ZZ_p}$ for $j \in \NN_0$.
Since $\zeta_{p,\,l,\,n}(0)$ corresponds to the related irregular pair of order $n$,
we then have
\[
   \zeta_{p,\,l,\,n}(0) \equiv \Delta_{(p,\,l)} \, s_{n+1} \pmod{p \ZZ_p}
\]
where $s_{n+1}$ is given by Definition \ref{def-char-padic}. A further consequence of
Proposition \ref{prop-irrpair-n}, extended to $s \in \ZZ_p$, is that
\begin{equation} \label{eqn-loc-zeta-padic-ordn-exp}
   \Dop \, \zeta_{p,\,l,\,n} (s) \equiv \Delta_{(p,\,l)} \pmod{p \ZZ_p}
\end{equation}
and furthermore
\[
   \zeta_{p,\,l,\,n} (s) \equiv \Delta_{(p,\,l)} \, ( s_{n+1} - s ) \pmod{p \ZZ_p} .
\]
By construction $\zeta_{p,\,l,\,n}$ is continuous; this is also a consequence of
\refeqn{eqn-loc-zeta-padic-ordn-diff-2}. Theorem \ref{theor-mahler} shows that
$\zeta_{p,\,l,\,n}$ has a Mahler series with the coefficients
$a_\nu = (-1)^\nu \, \Dop^\nu \zeta_{p,\,l,\,n}(0)$ for $\nu \geq 0$.
Using \refeqn{eqn-loc-zeta-padic-ordn-diff-1} we see that $a_\nu \in p^{n(\nu-1)}\ZZ_p$
for $\nu \geq 1$.
Thus, we can set $z_\nu = p^{-n(\nu-1)} a_\nu$ for $\nu \geq 1$ to obtain the proposed
series expansion above. From \refeqn{eqn-loc-zeta-padic-ordn-exp}
we deduce that $z_1 \equiv -\Delta_{(p,\,l)} \modp{p \ZZ_p}$.
\end{proof}

\begin{corl} \label{corl-zeta-chi-z}
Let $(p,l) \in \IRR_1$ with $\Delta_{(p,\,l)} \neq 0$.
Let
\[
   \zeta_{p,\,l,\,1}(s) = \zeta_{p,\,l,\,1}(0) +
     \sum_{\nu \geq 1} p^{\nu-1} z_\nu \binom{s}{\nu} ,
        \quad s \in \ZZ_p
\]
be the Mahler expansion of $\zeta_{p,\,l,\,1}$ given
by Proposition \ref{prop-zeta-padic-ordn-exp}. We have for $r \geq 1$ that
\[
   \zeta_{p,\,l}(s) \equiv \zeta_{p,\,l}(0) + \sum_{\nu = 1}^{r-1}
     p^\nu z_\nu \binom{s}{\nu} \pmod{p^r \ZZ_p} ,
        \quad s \in \ZZ_p .
\]
Special cases are given by
\[
   \zeta_{p,\,l}(0) = - \sum_{\nu \geq 1}
     p^\nu z_\nu \binom{\chi_{(p,\,l)}}{\nu}
\]
and for $r \geq 1$:
\[
   \zeta_{p,\,l}(0) \equiv - \sum_{\nu = 1}^{r-1}
     p^\nu z_\nu \binom{\chi_{(p,\,l)}}{\nu} \pmod{p^r \ZZ_p} .
\]

\end{corl}

\begin{proof}
We rewrite the Mahler expansion above using
$\zeta_{p,\,l,\,1}(s) = p^{-1} \zeta_{p,\,l}(s)$ which yields
\[
   \zeta_{p,\,l}(s) = \zeta_{p,\,l}(0) +
     \sum_{\nu \geq 1} p^\nu z_\nu \binom{s}{\nu} ,
        \quad s \in \ZZ_p .
\]
By the assumption $\Delta_{(p,\,l)} \neq 0$ we have the
zero $\chi_{(p,\,l)}$ of $\zeta_{p,\,l}$. The finite sums follow
easily $\mods{p^r \ZZ_p}$, since the coefficients $p^\nu z_\nu$
vanish for $\nu \geq r$.
\end{proof}

Hence, we can use Corollary \ref{corl-zeta-chi-z} to verify calculations
of the coefficients $z_\nu$ and of the zero $\chi_{(p,\,l)}$.
Now, we are almost ready to prove Theorem \ref{theor-zeta-padic-delta-kummer};
we first recall the following lemma from \cite[p.~227]{robert00}.

\begin{lemma} \label{lemma-binom-padic}
For $k \geq 1$ and $p^j \leq k < p^{j+1}$, we have
\[
   \left| \binom{s}{k} - \binom{t}{k} \right|_p \leq p^j \, |s-t|_p ,
     \quad s, t \in \ZZ_p .
\]
\end{lemma}

\begin{proofof}{Theorem~\ref{theor-zeta-padic-delta-kummer}}
Let $s,t \in \ZZ_p$ with $s \neq t$ and set $r = \ord_p ( s-t) \geq 0$.
Then we have by Kummer congruences that
\begin{equation} \label{eqn-loc-padic-dk-1}
   s \equiv t \pmod{p^{r}\ZZ_p} \quad \Longrightarrow \quad
     \zeta_{p,\,l} (s) \equiv \zeta_{p,\,l} (t) \pmod{p^{r+1}\ZZ_p} .
\end{equation}
We show that
\begin{equation} \label{eqn-loc-padic-dk-2}
   \frac{\zeta_{p,\,l} (s) - \zeta_{p,\,l} (t)}{p \, ( s - t )}
      \equiv - \Delta_{(p,\,l)} \pmod{p\ZZ_p}
         \quad \mbox{for} \quad s \neq t
\end{equation}
which also implies the converse to \refeqn{eqn-loc-padic-dk-1}.
Set $n = \min \, \{ \ord_p \zeta_{p,\,l} (s)$, $\ord_p \zeta_{p,\,l} (t) \}$
where $n \geq 1$. By Definition \ref{def-zeta-padic-ordn} we rewrite
$s = \psi_{n-1}(\chi_{(p,\,l)}) \, + p^{n-1} \, s'$
and $t = \psi_{n-1}(\chi_{(p,\,l)}) + p^{n-1} \, t'$ with some $s',t' \in \ZZ_p$.
Set $u=r+1-n \geq 0$. Note that $\ord_p (s'-t') = u$.
By the implication in \refeqn{eqn-loc-padic-dk-1} we then have
\[
   p^{-n} ( \zeta_{p,\,l} (s) - \zeta_{p,\,l} (t) )
     \equiv \zeta_{p,\,l,\,n} (s') - \zeta_{p,\,l,\,n} (t')
     \equiv p^u c \pmod{p^{u+1} \ZZ_p}
\]
with some integer $c$. On the other side, we can use the Mahler expansion
of $\zeta_{p,\,l,\,n}$ given by Proposition \ref{prop-zeta-padic-ordn-exp}.
We obtain
\begin{eqnarray*}
   \zeta_{p,\,l,\,n} (s') - \zeta_{p,\,l,\,n} (t') &\equiv&
      z_1 \, (s'-t') + p^n \, z_2 \left[ \binom{s'}{2} - \binom{t'}{2} \right] \\
      && \quad  + \, p^{2n} \, z_3 \left[ \binom{s'}{3} - \binom{t'}{3} \right] + \ldots
        \quad \pmod{p^{u+1} \ZZ_p} .
\end{eqnarray*}
Since $\ord_p (s'-t') = u$, Lemma \ref{lemma-binom-padic} shows that
\[
   \ord_p \left[ \binom{s'}{k} - \binom{t'}{k} \right] \geq u - \lfloor \log_p k \rfloor
\]
where $\log_p$ is the real valued logarithm with base $p$.
Since we have $p \geq 5$ we obtain the estimate
\[
   \ord_p \left( p^{n(k-1)} z_k \left[ \binom{s'}{k} - \binom{t'}{k} \right] \right)
     \geq u + n(k-1) - \lfloor \log_p k \rfloor \geq u+1
\]
for all $k \geq 2$. Hence, all terms vanish $\mods{p^{u+1}}$ except for $k=1$:
\[
   \zeta_{p,\,l,\,n} (s') - \zeta_{p,\,l,\,n} (t') \equiv z_1 \, (s'-t') \pmod{p^{u+1} \ZZ_p} .
\]
By Proposition \ref{prop-zeta-padic-ordn-exp} and $\ord_p (s'-t') = u$ we finally get
\[
   z_1 \, (s'-t') \equiv - \Delta_{(p,\,l)} (s'-t') \equiv p^u c \pmod{p^{u+1} \ZZ_p}
\]
which shows that $c \not\equiv 0 \modp{p}$.
Since $p^{-n+1} \, (s-t) = s'-t'$, this deduces \refeqn{eqn-loc-padic-dk-2}.
\smallskip

Now, we have to determine the derivative of $\zeta_{p,\,l}$.
Taking any sequence $(t_\nu)_{\nu \geq 1}$ with $t_\nu \neq s$ for all $\nu \geq 1$
and $\lim_{\nu \rightarrow \infty} t_\nu = s$,
\refeqn{eqn-loc-padic-dk-2} yields the derivative with
\[
   \zeta'_{p,\,l} (s) \equiv - p \, \Delta_{(p,\,l)} \pmod{p^2\ZZ_p} .  \vspace*{-3ex}
\]
\end{proofof}

\begin{remark}
In general the converse of the Kummer congruences does not hold.
The first nontrivial counterexample is given by $p=13$ and
$B_{16}/16 - B_4/4 = -7\cdot13^2/2720$.
The only case where the divided Bernoulli numbers are equal
is $B_{14}/14 - B_2/2 = 0$.
\end{remark}

At the end of this section we revisit the $p$-adic zeta functions $\zeta_{p,\,s_1}$
in the case $s_1 \neq 0$ as introduced by Definition \ref{def-zeta-padic}.
We can transfer some results to these functions.
Moreover, we can state a formula equivalent to the Kummer congruences,
but which involves values of $p$-adic zeta functions
at smallest possible argument values.

\begin{prop} \label{prop-zeta-padic-mahler}
Let $p$ be a prime with $p \geq 5$ and
$s_1 \in \{ 2,4,\ldots,p-3 \}$ be fixed. Then
\begin{equation} \label{eqn-loc-zeta-padic-mahler-1}
   \Dop^{r} \zeta_{p,\,s_1} (s) \equiv 0 \pmod{p^r \ZZ_p} ,
     \quad s \in \ZZ_p .
\end{equation}
There exists the Mahler expansion
\begin{equation} \label{eqn-loc-zeta-padic-mahler-2}
   \zeta_{p,\,s_1}(s) = \zeta_{p,\,s_1}(0) +
     \sum_{\nu \geq 1} p^{\nu} z_\nu \binom{s}{\nu} ,
        \quad s \in \ZZ_p
\end{equation}
with $z_\nu \in \ZZ_p$
where $z_\nu = (-1)^\nu p^{-\nu} \, \Dop^\nu \zeta_{p,\,s_1}(0)$.
\end{prop}

\begin{proof}
In analogy to Propositions \ref{prop-zeta-padic-ordn-diff} and
\ref{prop-zeta-padic-ordn-exp} in the case $n=1$,
Congruence \refeqn{eqn-loc-zeta-padic-mahler-1} is a
consequence of Theorem \ref{theor-kummer-congr-general}.
Since $\zeta_{p,\,s_1}$ is a continuous function on $\ZZ_p$, we obtain
a Mahler expansion where the proposed coefficients follow by
\refeqn{eqn-loc-zeta-padic-mahler-1}.
\end{proof}

\begin{theorem} \label{theor-kummer-congr-zeta}
Let $p$ be a prime with $p \geq 5$ and
$s_1 \in \{ 2,4,\ldots,p-3 \}$ be fixed. For $r \geq 0$ we have
\[
   \zeta_{p,\,s_1}(s) \equiv \sum_{k = 0}^{r} \zeta_{p,\,s_1}(k)
     \,\, T_{r,k}(s) \pmod{p^{r+1} \ZZ_p} ,
        \quad s \in \ZZ_p
\]
where the polynomials $T_{r,k}: \ZZ_p \to \ZZ_p$ with
$T_{r,k} \in \QQ[x]$ and $\deg T_{r,k} = r$ are given by
\[
   T_{r,k}(x) = \sum_{j=k}^{r} (-1)^{j+k} \binom{j}{k} \binom{x}{j} .
\]
\end{theorem}

\begin{proof}
We rewrite the Mahler expansion of $\zeta_{p,\,s_1}$ given in
\refeqn{eqn-loc-zeta-padic-mahler-2}.
Let $s \in \ZZ_p$. For $r \geq 0$ we obtain the finite expansion
\[
   \zeta_{p,\,s_1}(s) \equiv \zeta_{p,\,s_1}(0) +
       \sum_{\nu=1}^r p^{\nu} z_\nu \binom{s}{\nu}
     \equiv \sum_{j=0}^r (-1)^j \, \Dop^j \zeta_{p,\,s_1}(0) \binom{s}{j}
     \pmod{p^{r+1} \ZZ_p} .
\]
By Definition \ref{def-diff-oper} we have
\[
   \Dop^j \zeta_{p,\,s_1}(0)
     = \sum_{k=0}^j \binom{j}{k} (-1)^{k} \zeta_{p,\,s_1}(k)
     = \sum_{k=0}^r \binom{j}{k} (-1)^{k} \zeta_{p,\,s_1}(k) .
\]
Reordering the finite sums and terms yields the definition of $T_{r,k}$ above.
Since $x \mapsto \binom{x}{j}$ which is a polynomial of degree $j$
defines a function on $\ZZ_p$, the polynomials $T_{r,k}$ are functions on $\ZZ_p$.
The coefficients of $T_{r,k}$ are rational, e.g., rewrite
$T_{r,r}(x) = \binom{x}{r}$ as a polynomial in $x$.
We deduce that $\deg T_{r,k} = r$, because
the term $\binom{x}{r}$ which occurs only once gives the maximal degree.
\end{proof}

\begin{corl} \label{corl-kummer-congr-zeta}
Let $n$ be an even positive integer and $p$ be a prime where $p-1 \notdiv n$.
Define the integer $s_1$ by $s_1 \equiv n \modp{p-1}$ with $0 < s_1 < p-1$.
Set $s = (n-s_1)/(p-1)$. For $r \geq 1$ we have
\[
   (1-p^{n-1}) \BN(n) \equiv \sum_{k = 0}^{r-1}
      (1-p^{s_1+k(p-1)-1}) \BN(s_1+k(p-1))
        \,\, T_{r-1,k}(s) \!\! \pmod{p^r \ZZ_p}
\]
with $T_{r,k}$ as defined above.
\end{corl}

\begin{proof}
This is a reformulation of Theorem \ref{theor-kummer-congr-zeta} using
Definition \ref{def-zeta-padic} where $s_1 + s(p-1) = n$.
\end{proof}

\begin{remark}
The case $r=1$ of Corollary \ref{corl-kummer-congr-zeta}
reduces to a special case of the Kummer congruences \refeqn{eqn-kummer-congr} for $r=1$.
We explicitly give the cases $r=2,3,4$ of Corollary \ref{corl-kummer-congr-zeta}.
Note that $s_1 \geq 2$ and $p \geq 5$.
Case $r=2$:
\[
   (1-p^{n-1}) \BN(n) \equiv
     - (s-1) \, (1-p^{s_1-1}) \BN(s_1)
     + s \, \BN(s_1+p-1) \pmod{p^2 \ZZ_p} .
\]
Case $r=3$:
\begin{eqnarray*}
   (1-p^{n-1}) \BN(n) &\equiv&
        + \, \frac12 (s^2-3s+2) \, (1-p^{s_1-1}) \BN(s_1) \\
     && - \, s(s-2) \, \BN(s_1+p-1) \\
     && + \, \frac12 s(s-1) \, \BN(s_1+2(p-1)) \qquad
        \pmod{p^3 \ZZ_p} .
\end{eqnarray*}
Case $r=4$:
\begin{eqnarray*}
   (1-p^{n-1}) \BN(n) &\equiv&
        - \frac16 (s^3-6s^2+11s-6) \, (1-p^{s_1-1}) \BN(s_1) \\
     && + \, \frac12 s(s^2-5s+6) \, \BN(s_1+p-1) \\
     && - \, \frac12 s(s^2-4s+3) \, \BN(s_1+2(p-1)) \\
     && + \, \frac16 s(s-1)(s-2) \, \BN(s_1+3(p-1)) \qquad
        \pmod{p^4 \ZZ_p} .
\end{eqnarray*}
\end{remark}

\section{Algorithms}
\setcounter{equation}{0}

Here we will give some algorithms for calculating irregular pairs of higher order
assuming we are in the nonsingular case $\Delta_{(p,\,l)} \neq 0$.
As a result of Theorem \ref{theor-delta-irr} and Theorem \ref{theor-zeta-padic},
one has first to calculate $\Delta_{(p,\,l)}$, then the irregular pair of order $n$
resp.\ the corresponding divided Bernoulli number
provides the next related irregular pair of order $n+1$.
This is not practicable for higher orders say $n>3$.
Proposition \ref{prop-irrpair-rn} shows a more effective way to determine
related irregular pairs of higher order. Starting from an irregular pair
$(p,l) \in \IRR_n$ with $n \geq 1$ and requiring that $l > (r-1)n$
with some $r \geq 2$
we can obtain a related irregular pair $(p,l') \in \IRR_{rn}$.
If the corresponding sequence
$( \alpha_j )_{j \geq 0}$ is equidistant $\mods{p^{(r-1)n}}$, then one
can easily apply this proposition. If not, one has to calculate successively
all elements $\alpha_j$ in order to find $\alpha_s \equiv 0 \modp{p^{(r-1)n}}$
where $0 \leq s < p^{(r-1)n}$ which exists uniquely by the assumption $\Delta_{(p,\,l)} \neq 0$.
To shorten these calculations, this search can be accomplished step by step,
moving each time from a sequence $( \alpha_{j,k} )_{j \geq 0}$ to
a sequence $( \alpha_{j,k+1} )_{j \geq 0}$ which are assigned to the irregular pair
of order $k$ resp.\ $k+1$.

\begin{prop} \label{prop-calc-irrpair}
Let $(p,l) \in \IRR_n$, $n \geq 1$ with $\Delta_{(p,\,l)} \neq 0$.
Let $r$, $u$ be positive integers with
$r>1$ and $u=(r-1)n$. Assume that $l>u$. Let the elements
\[
    \alpha_{j,0} \equiv p^{-n} \BN(l + j \eulerphi(p^n))
      \pmod{p^u} , \quad  j = 0, \ldots, r-1
\]
be given. For each step $k=0,\ldots,u-1$ proceed as follows:
\smallskip

\noindent The elements $\alpha_{j,k}$ with $j=0,\ldots,rp-1$ have to be calculated
successively by
\[
    \alpha_{j+r,k} \equiv (-1)^{r+1} \sum_{\nu=0}^{r-1} \binom{r}{\nu} (-1)^{\nu}
      \alpha_{j+\nu,k} \pmod{p^{u-k}} .
\]
Set $s_k \equiv -\alpha_{0,k} \, \Delta_{(p,\,l)}^{-1} \modp{p}$ with $0 \leq s_k < p$.
The elements $\alpha_{j,k}$ which are divisible by $p$ are given by
\[
    \alpha_{s_k+\mu p,k} \equiv 0 \pmod{p} , \quad \mu=0,\ldots,r-1 .
\]
For $k < u-1$ set
\[
   \alpha_{j,k+1} = \alpha_{s_k+jp,k} / p , \quad j=0,\ldots,r-1
\]
and go to the next step $k+1$, otherwise stop.
Let $(p,t_1, \ldots, t_n) \in \IRP_n$ be the associated irregular pair with $(p,l)$, then
\[
   (p,t_1, \ldots, t_n, s_0, \ldots, s_{u-1}) \in \IRP_{rn}
\]
is the only related irregular pair of order $rn$.
\end{prop}

\begin{proof}
Proposition \ref{prop-irrpair-rn} shows that
\begin{equation} \label{eqn-loc-irrpair-calc-1}
   \sum_{\nu=0}^r \binom{r}{\nu} (-1)^{\nu}
      \alpha_{j+\nu,0} \equiv 0  \pmod{p^{u}} .
\end{equation}
All elements of the sequence $(\alpha_{j,0})_{j \geq 0}$ can be calculated
successively induced by the first elements $\alpha_{j,0}$ with $j = 0, \ldots, r-1$.
Using Corollary \ref{corl-kummer-congr-irrpair} with $\omega = p^k\eulerphi(p^n)$, $0 \leq k < u$,
then \refeqn{eqn-loc-irrpair-calc-1} becomes
\begin{equation} \label{eqn-loc-irrpair-calc-2}
   \sum_{\nu=0}^r \binom{r}{\nu} (-1)^{\nu}
      \alpha_{j+\nu p^k,0} \equiv 0  \pmod{p^{u}} ,
\end{equation}
whereby the sequence $(\alpha_{j+\mu p^k,0})_{\mu \geq 0}$ can also
be calculated successively. Note that now $j$ is fixed and $\mu$ runs.
The sequences $(\alpha_{j,k})_{j \geq 0}$ which we will consider are
subsequences of $(\alpha_{j,0})_{j\geq0}$ in a suitable manner.
Essentially, these sequences are given by \refeqn{eqn-loc-irrpair-calc-2}.
The existence of these sequences and that they correspond to the related irregular pair
of order $n+k$ will be shown by induction on $k$ for $k=0,\ldots,u-1$.
Set $l_n=l$. By Proposition \ref{prop-irrpair-n} and Theorem \ref{theor-delta-irr}
there exist certain integers $s_k$ and related irregular pairs of higher order
for $k=0,\ldots,u-1$ where
\[
    (p,l_{n+k+1}) \in \IRR_{n+k+1} ,
      \quad l_{n+k+1} = l_{n+k} + s_k \, \eulerphi(p^{n+k})
      \quad \mbox{with} \quad 0 \leq s_k < p .
\]
Basis of induction $k=0$: The sequence $(\alpha_{j,0})_{j\geq0}$ is given
by \refeqn{eqn-loc-irrpair-calc-1} and we have
\[
    \alpha_{j,0} \equiv p^{-n} \, \BN(l_n + j \eulerphi(p^n))
      \pmod{p^u} .
\]
Inductive step $k \mapsto k+1$: Assume true for $k$ prove for $k+1$.
The elements $\alpha_{j,k}$ with $j=0,\ldots,r-1$ are given and the
following elements are calculated by
\begin{equation} \label{eqn-loc-irrpair-calc-3}
    \alpha_{j+r,k} \equiv (-1)^{r+1} \sum_{\nu=0}^{r-1} \binom{r}{\nu} (-1)^{\nu}
      \alpha_{j+\nu,k} \pmod{p^{u-k}}
\end{equation}
up to index $j=rp-1$. Proposition \ref{prop-irrpair-n} provides
\begin{equation} \label{eqn-loc-irrpair-calc-4}
   s_k \equiv -\alpha_{0,k} \Delta_{(p,\,l)}^{-1} \pmod{p}
     \quad \mbox{with} \quad 0 \leq s_k < p .
\end{equation}
In the case $k < u-1$, Lemma \ref{lemma-seq-a2} additionally ensures
that only $\alpha_{s_k+jp,k} \equiv 0 \modp{p}$ for all $j$.
Thus, we can define a new sequence by
\begin{eqnarray} \nonumber
    \alpha_{j,k+1} &\equiv& \alpha_{s_k+jp,k} / p \\
  \nonumber
     &\equiv& p^{-(n+k+1)} \, \BN(l_{n+k} + (s_k+jp) \eulerphi(p^{n+k})) \\
  \label{eqn-loc-irrpair-calc-5}
     &\equiv& p^{-(n+k+1)} \, \BN(l_{n+k+1} + j \eulerphi(p^{n+k+1}))
       \qquad \pmod{p^{u-(k+1)}}
\end{eqnarray}
for $j=0,\ldots,r-1$. By definition $(p\,\alpha_{j,k+1})_{j=0,\ldots,r-1}$ is
a subsequence of $(\alpha_{j,k})_{j\geq0}$. Inductively
$(p^{k+1}\,\alpha_{j,k+1})_{j=0,\ldots,r-1}$ is a subsequence
of $(\alpha_{j,0})_{j\geq0}$ and satisfies \refeqn{eqn-loc-irrpair-calc-2}
in a suitable manner and therefore
also \refeqn{eqn-loc-irrpair-calc-3} considering case $k+1$.
On the other side, Congruence \refeqn{eqn-loc-irrpair-calc-5}
shows that the new sequence also corresponds to the related
irregular pair of order $n+k+1$.

Let $(p,t_1, \ldots, t_n) \in \IRP_n$ be the associated irregular pair with $(p,l)$.
Congruence \refeqn{eqn-loc-irrpair-calc-4} provides a unique integer $s_k$
for each step. Thus,
$(p,t_1, \ldots, t_n, s_0, \ldots, s_{u-1}) \in \IRP_{rn}$
is the only related irregular pair of order $rn$.
\end{proof}

\begin{remark}
Unfortunately, the previous proposition has the restriction that
for an irregular pair $(p,l) \in \IRR_n$ of order $n$ and
parameter $r$ the following must hold that
\[
    l > (r-1) n .
\]
Consider $(691,12) \in \IRR_1$. In this case one only could calculate related irregular pairs
up to order 12. However, this restriction can be removed by shifting the index of the
starting sequence $(\alpha_{j,0})_{j \geq 0}$. Then shifting $j \mapsto j+t$ yields
\[
   l + t \, \eulerphi(p^n) > (r-1) n
\]
enabling one to choose a greater value of $r$. In general, one has to proceed in
the following way. Moving from a sequence $(\alpha_{j,k})_{j \geq 0}$ to
$(\alpha_{j,k+1})_{j \geq 0}$ one has to determine elements
$\alpha_{j,k} \equiv 0 \modp{p}$. If one starts with a shifted sequence
$(\alpha'_{j,k})_{j \geq 0} = (\alpha_{j,k})_{j \geq t}$ having elements
$\alpha_{j,k} \equiv 0 \modp{p}$ with $0 \leq j < t$, then
the calculated sequence $(\alpha'_{j,k+1})_{j \geq 0}$
is also shifted in the index compared to $(\alpha_{j,k+1})_{j \geq 0}$.
This can be easily corrected by comparing the sequences
resp.\ the resulting integers $s_k$ with
unshifted sequences calculated with a lower $r' < r$.
In this case determining the integer $s_k$ is better done by searching
$\alpha_{j,k} \equiv 0 \modp{p}$ rather than calculating via \refeqn{eqn-loc-irrpair-calc-4}.
\end{remark}

The main result can be stated as follows: irregular pairs of higher order can be determined with
little effort by calculating a small number of divided Bernoulli numbers with small indices.
We shall give another algorithm in terms of the $p$-adic zeta function $\zeta_{p,\,l}$
which produces a truncated $p$-adic expansion of $\chi_{(p,\,l)}$.

\begin{prop} \label{prop-calc-chi}
Let $(p,l) \in \IRR_1$ with $\Delta_{(p,\,l)} \neq 0$.
Let $n$ be a positive integer. Initially calculate the values
\[
   \zeta_{p,\,l,\,1}(k) \equiv p^{-1} \zeta_{p,\,l}(k) \pmod{p^n \ZZ_p}
\]
for $k=0,\ldots,n$ and
\[
   \Delta_{(p,\,l)} \equiv \zeta_{p,\,l,\,1}(0) - \zeta_{p,\,l,\,1}(1) \pmod{p \ZZ_p}
\]
where $0 < \Delta_{(p,\,l)} < p$. Set $t_1=0$.
For each step $r=1,\ldots,n$ proceed as follows:
\smallskip

\noindent Calculate
\[
   \xi_r \equiv \sum_{k=0}^r \zeta_{p,\,l,\,1}(k) \, T_{r,k}(t_r) \pmod{p^r \ZZ_p}
\]
with the polynomials $T_{r,k}$ as defined in Theorem \ref{theor-kummer-congr-zeta}.
Then $\xi_r \in p^{r-1} \ZZ_p$. Set
\[
   s_{r+1} \equiv \Delta_{(p,\,l)}^{-1} \, p^{1-r} \, \xi_r \pmod{p \ZZ_p}
\]
where $0 \leq s_{r+1} < p$. Set $t_{r+1} = t_r + s_{r+1} \, p^{r-1}$ and
go to the next step while $r < n$.
Finally $t_{n+1} \equiv \chi_{(p,\,l)} \modp{p^n \ZZ_p}$.
\end{prop}

\begin{proof}
By Definition \ref{def-zeta-padic-ordn} we have
$\zeta_{p,\,l,\,1}(k) = p^{-1} \zeta_{p,\,l}(k)$. The value
$\Delta_{(p,\,l)}$ is given by its definition.
Define $t'_r = \psi_{r-1}( \chi_{(p,\,l)} )$ for any $r \geq 1$
where the expansion of the zero of $\zeta_{p,\,l}$ is given by
\[
   \chi_{(p,\,l)} = s_2 + s_3 \, p + \ldots + s_{r+1} \, p^{r-1} + \ldots \,.
\]
For now, let $r \in \{1,\ldots,n\}$ be fixed.
Theorem \ref{theor-kummer-congr-zeta} provides
\[
   \zeta_{p,\,l,\,1}(s) \equiv \sum_{k = 0}^{r} \zeta_{p,\,l,\,1}(k)
     \,\, T_{r,k}(s) \pmod{p^{r} \ZZ_p}
\]
for $s \in \ZZ_p$.
From Definition \ref{def-zeta-padic-ordn} we have
\[
   \zeta_{p,\,l,\,r}(0) = p^{-r} \zeta_{p,\,l}( \psi_{r-1}( \chi_{(p,\,l)} ) )
     = p^{1-r} \zeta_{p,\,l,\,1}( t'_r ) .
\]
Proposition \ref{prop-zeta-padic-ordn-exp} shows that
\[
   \zeta_{p,\,l,\,r}(0) \equiv \Delta_{(p,\,l)} \, s_{r+1} \pmod{p \ZZ_p} .
\]
By construction we have $t'_{r+1} = t'_r + s_{r+1} \, p^{r-1}$.
Since $t'_1 = t_1 = 0$ we deduce by induction on $r$ that
$t'_r = t_r$ and $\xi_r \equiv \zeta_{p,\,l,\,1}( t'_r ) \modp{p^r \ZZ_p}$
for $1 \leq r \leq n$.
This produces $t_{n+1} = \psi_n( \chi_{(p,\,l)} )$ which
equals to $t_{n+1} \equiv \chi_{(p,\,l)} \modp{p^n \ZZ_p}$.
Since $1 \leq r \leq n$, we need the values
$\zeta_{p,\,l,\,1}(k) \modp{p^n \ZZ_p}$ for $k=0,\ldots,n$.
\end{proof}

\begin{prop} \label{prop-bn-repr}
Let $n$ be an even positive integer, then
\[
   B_n = (-1)^{\frac{n}{2}-1} \,
    \prod\limits_{p-1 \notdiv n} p^{\tau(p,n)+\ord_p n} \,\,\Big/\,\,
     \prod\limits_{p-1 \pdiv n} p
\]
where
\[
   \tau(p,n) := \sum_{\nu=1}^\infty \#( \, \IRR_\nu \cap
     \{ ( p, n \modres \eulerphi(p^\nu) ) \} \, ) .
\]
Here, as in Definition \ref{def-ordn-irrpair}, $x \modres y$ denotes the least nonnegative
residue of $x$ modulo $y$.
\end{prop}

\begin{proof}
The trivial parts of the products above are given
by \refeqn{eqn-clausen-von-staudt}, \refeqn{eqn-adams-triv-div}, and the sign.
Thus, the product $\prod_{p-1 \notdiv n} p^{\tau(p,n)}$ remains.
From Definition \ref{def-ordn-irrpair} and Remark \ref{rem-ordn-irrpair},
the function $\tau(p,n)$ follows by applying the maps $\lambda_\nu$ resp.\
Kummer congruences which results in a chain of related irregular pairs of descending order
similar to \refeqn{eqn-chain-irrpair}.
\end{proof}

The previous proposition gives an unconditional representation of the Bernoulli numbers
by means of the sets $\IRR_\nu$. Theorem \ref{theor-zeta-prod-psi} also gives a representation
by zeros $\chi_{(p,\,l)}$ assuming the $\Delta$-Conjecture.
Of course, the problem of determining the occurring irregular prime factors remains open.
On the other side, for instance, if one has calculated the first irregular pairs
of order 10 for the first irregular primes
$p_1, \ldots, p_r$ like Table \ref{tab-irrpair-ord-10}, then one can specify
\textsl{ad hoc} all irregular prime powers $p_\nu^{e_\nu}$ with $p_\nu \leq p_r$ of $B_n$
resp.\ $\zeta(1-n)$ up to index $n=4 \!\cdot\! 10^{15}$. Note that this lower bound is
here determined by the first irregular prime 37 and order 10.
\smallskip

Define for positive integers $n$ and $m$ the summation function of consecutive integer powers by
$S_n(m) = \sum_{\nu=0}^{m-1} \nu^n$. Many congruences concerning the function $S_n$
are naturally related to the Bernoulli numbers.

\begin{prop} \label{prop-delta-sn}
Let $(p,l) \in \IRR_1$. Assume that $n \equiv l \modp{p-1}$ where $n > 0$, then
\begin{equation} \label{eqn-delta-sn-1}
   \frac{B_n}{n} \equiv \frac{S_n(p)}{n \, p} \pmod{p^2} .
\end{equation}
Moreover
\begin{equation} \label{eqn-delta-sn-2}
   \Delta_{(p,\,l)} \equiv p^{-2} \left(
     \frac{S_{l+p-1}(p)}{l-1} - \frac{S_{l}(p)}{l} \right) \pmod{p}
\end{equation}
with $0 \leq \Delta_{(p,\,l)} < p$.
\end{prop}

\begin{proof}
Let $n \equiv l \modp{p-1}$.
The well-known formula of $S_n$, see \cite[p.~234]{ireland90}, is given by
\begin{equation} \label{eqn-loc-delta-sn}
   \frac{S_n(p)}{n \, p} = \frac{B_n}{n} + \binom{n-1}{1} B_{n-2} \frac{p^{2}}{2 \cdot 3} +
      \sum_{k=3}^n \binom{n-1}{k-1} B_{n-k} \frac{p^{k}}{k(k+1)}
\end{equation}
where the equation is divided by $n$ and $p$. Note that $p \geq 37$ and $n \geq l \geq 12$, because
37 is the first irregular prime and $B_{12}/12$ is the first divided Bernoulli number which has
a numerator greater one. Now, the properties of \refeqn{eqn-clausen-von-staudt} and
\refeqn{eqn-adams-triv-div} provide that $B_n/n$ and $B_{n-2}$ are $p$--integral.
For the other terms with $B_{n-k} \neq 0$ it follows that $p B_{n-k}$ is $p$--integral and
$\ord_p \, (p^{k-1}/(k(k+1))) \geq 2$ by a standard counting argument. Therefore,
Equation \refeqn{eqn-loc-delta-sn} is $p$-integral and holds $\mods{p^2}$ whereas
all terms of the right side vanish except for $B_n/n$.
This gives Congruence \refeqn{eqn-delta-sn-1}.
From Definition \ref{def-delta} we have
\[
   p \, \Delta_{(p,\,l)} \, \equiv \, \BN(l + p-1) - \BN(l) \pmod{p^2} .
\]
We can apply \refeqn{eqn-delta-sn-1} in the congruence above.
Reducing a $p$-power and considering that $l+p-1 \equiv l-1 \not\equiv 0 \modp{p}$
finally yields \refeqn{eqn-delta-sn-2}.
\end{proof}

Looking at each line of Table \ref{tab-irrpair-ord-10},
the product of the first three entries $\Delta_{(p,\,l)}$, $s_1$, and $s_2$ are connected with
the function $S_n$. Thus, one can easily verify these values.

\begin{prop} \label{prop-delta-s1-s2}
Let $(p,l) \in \IRR_1$ with $\Delta_{(p,\,l)} \neq 0$.
Let $(p,s_1,s_2) \in \IRP_2$ be the related
irregular pair of order two with $l=s_1$. Then
\[
   \Delta_{(p,\,l)} \, s_1 \, s_2 \equiv - p^{-2} \, S_l(p) \pmod{p} .
\]
\end{prop}

\begin{proof}
By Proposition \ref{prop-irrpair-n} we have
\[
   s_2 \equiv - p^{-1} \frac{B_l}{l} \Delta_{(p,\,l)}^{-1} \pmod{p} .
\]
We then obtain by Proposition \ref{prop-delta-sn} that
\[
   p \, s_2 \, \Delta_{(p,\,l)} \equiv - \frac{B_l}{l}
     \equiv - \frac{S_l(p)}{l \, p} \pmod{p^2}
\]
which deduces the result since $s_1 = l < p$.
\end{proof}

Now, we shall give some reasons why a prediction or description of
the occurrence of irregular prime factors of
Bernoulli numbers seems to be impossible in general. For example, we have with an
extremely small index $n=42$ that
\[
    B_{42} = \frac{1520097643918070802691}{1806}
\]
observing that the numerator is a large irregular prime! As mentioned in Section \ref{sect-main},
Bernoulli numbers $B_n$ with index up to $n=152$ have large irregular prime factors
with 30 up to 100 digits. This is even now the greatest mystery of the Bernoulli numbers!

The connection with the Riemann zeta function $\zeta(s)$ via \refeqn{eqn-zeta-n} leads
to methods of calculating Bernoulli numbers directly in a fast and effective way,
see \cite[Section 2.7]{kellner02}, noting that
the main part of the calculation can be done using integers only.
Let $n$ be an even positive integer and $|B_n|=U_n/V_n$ with
$(U_n,V_n)=1$ then \refeqn{eqn-zeta-n} reads
\[
   U_n = \tau_n \, \zeta(n) , \qquad \tau_n = 2 V_n \frac{n!}{(2\pi)^n} ,
     \qquad V_n = \prod_{p-1 \pdiv n} p
\]
with $V_n$ given by \refeqn{eqn-clausen-von-staudt}. Since $\zeta(n) \rightarrow 1$
for $n \rightarrow \infty$, $\tau_n$ is a first approximation of the numerator $U_n$.
Considering the decimal digit representation of $U_n$ and $\tau_n$,
about $n/3$ digits of the most significant decimal digits of
$U_n$ and $\tau_n$ are equal, see \cite[Satz 2.7.9, p.~75]{kellner02} for a
more precise statement and formula.
Visiting $B_{42}$ again, we observe 12 identical digits:
\begin{eqnarray*}
 U_{42} &=& \textbf{1\,520\,097\,643\,91}8\,070\,802\,691 , \\
 \tau_{42} = 2 \cdot 1806 \cdot 42! / (2\pi)^{42} &\approx&
   \textbf{1\,520\,097\,643\,91}7\,725\,172\,488.7773 \,.
\end{eqnarray*}

How can we interpret this result? A part of the most significant digits
of the numerator of the Bernoulli number $B_n$ is determined in a certain way
by all primes $p \le n+1$ and the reciprocal $n$-th power of $\pi$.
That the digits of $\pi$ are involved in the numerators of the Bernoulli numbers
is quite remarkable.

\section{Connections with Iwasawa theory}
\setcounter{equation}{0}
\label{sect-iwasawa}

The $\Delta$-Conjecture is directly connected with Iwasawa theory
of cyclotomic fields over $\QQ$. Let $p$ be an odd prime and
$\mu_{p^n}$ be the set of $p^n$-th roots of unity where $n$ is a positive
integer. For the cyclotomic field $\QQ(\mu_{p^n})$ let $\QQ(\mu_{p^n})^+$ be
its maximal real subfield. The class number
$h_p = h(\QQ(\mu_p))$ can be factored by $h_p = h_p^- \, h_p^+$
where $h_p^+ = h(\QQ(\mu_{p})^+)$ and $h_p^-$ is
the relative class number introduced by Kummer.
Define
\[
    B_{1,\omega^{m}} = \frac{1}{p} \sum_{a=1}^{p-1} a \, \omega^{m}(a)
\]
as the generalized Bernoulli number assigned to the Teichm\"uller character $\omega$.
This character is defined by $\omega : \, \ZZ_p^* \rightarrow \ZZ_p^*$
and $\omega(a) \equiv a \modp{p\ZZ_p}$ for $a \in \ZZ_p^*$ giving the
$(p-1)$-th roots of unity in $\QQ_p$.
We have for even $m > 0$ and $p-1 \notdiv m$ the following relation,
see \cite[Corollary 5.15, p. 61]{washington97}:
\[
    B_{1,\omega^{m-1}} \equiv \BN(m) \pmod{p\ZZ_p} .
\]

For the detailed theory, especially of Iwasawa invariants and cyclotomic
$\ZZ_p$-exten\-sions, see Washington \cite{washington97} and
Greenberg \cite{greenberg01}.
The results of Iwasawa, Ferrero and Washington, Vandiver and Kummer
provide the following theorem, see \cite[Cor.~10.17, p.~202]{washington97}.

\begin{theorem} \label{theor-iwasawa}
Let $p$ be an irregular prime. Assume the following conditions
for all irregular pairs $(p,l)$:
\begin{enumerate}
\item The conjecture of Kummer--Vandiver holds: $p \notdiv h_p^+$,
\item The Kummer congruence does not hold $\mods{p^2}$:
  $\BN(l+p-1) \not\equiv \BN(l) \modp{p^2}$,
\item The generalized Bernoulli number is not divisible by $p^2$:
  $B_{1,\omega^{l-1}} \not\equiv 0 \modp{p^2\ZZ_p}$.
\end{enumerate}
If these are satisfied, then
\[
   \ord_p h( \QQ( \mu_{p^n} ) ) \, = \, i(p) \, n \quad \mbox{for all} \quad n \geq 1 .
\]
\end{theorem}

All conditions of the theorem above hold for all irregular primes $p < 12\,000\,000$
as verified in \cite{irrprime12M}.
In the case of a regular prime $p$ the formula of the theorem above is also valid,
because then we have $i(p)=0$ and
\[
    p \notdiv h_p = h(\QQ(\mu_p)) \quad \Longleftrightarrow \quad
      p \notdiv h(\QQ(\mu_{p^n})) \quad \mbox{for all} \quad n \geq 1
\]
as follows, e.g., from Iwasawa theory. To get another point of view we can exchange
two conditions of the previous theorem by our results.
Conditions equivalent to those of Theorem \ref{theor-iwasawa} are as follows:

\begin{enumerate}
\item[(2')] The $\Delta$-Conjecture holds: $\Delta_{(p,\,l)} \neq 0$,
\item[(3')] A special irregular pair of order two does not exist: $(p,l,l-1) \notin \IRP_2$.
\end{enumerate}

Now, the $\Delta$-Conjecture with its consequences gives a significant reason to believe
that the condition (2) resp.\ (2') may hold in general.
We still have to show that the condition (3') is equivalent to the condition (3).

\begin{prop} \label{prop-bw-bl-2}
Let $(p,l) \in \IRR_1$. Then
\[
   B_{1,\omega^{l-1}} \equiv B_{l+(p-1)(l-1)} \pmod{p^2\ZZ_p}
\]
and
\[
   B_{1,\omega^{l-1}} \equiv 0 \pmod{p^2\ZZ_p} \quad
     \Longleftrightarrow \quad (p,l,l-1) \in \IRP_2 .
\]
\end{prop}

To prove this result, we first need some properties of the Teichm\"uller character $\omega$.
Since $\omega(a)$ is defined by $\omega(a) = \lim_{n \to \infty} a^{p^n}$ in $\ZZ_p$,
the following lemma is easily derived.

\begin{lemma} \label{lemma-teich}
Let $a$, $p$ be integers with $p$ an odd prime and $0 < a < p$. Then
\[
   \omega(a) \,\equiv\, a^p \pmod{p^2\ZZ_p} \qquad \mbox{and} \qquad
     \omega(a) \,\equiv\, a^p + p \, (a^p-a) \pmod{p^3\ZZ_p} .
\]
\end{lemma}
\medskip

\begin{proofof}{Proposition~\ref{prop-bw-bl-2}}
Since $\omega^{l-1}(a) = \omega(a^{l-1})$, we have
\[
    p \, B_{1,\omega^{l-1}} = \sum_{a=1}^{p-1} a \, \omega^{l-1}(a)
      = \sum_{a=1}^{p-1} a \, \omega(a^{l-1}) .
\]
Using Lemma \ref{lemma-teich} we obtain
\[
    a \, \omega(a^{l-1}) \equiv a^{p(l-1)+1} + p\,( a^{p(l-1)+1} - a^l ) \pmod{p^3\ZZ_p} .
\]
From the definition of $S_n$ and Proposition \ref{prop-delta-sn} we have
\[
    p B_{1,\omega^{l-1}} \equiv S_{p(l-1)+1}(p) + p S_{p(l-1)+1}(p) - p S_l(p)
      \equiv p B_{l + (p-1)(l-1)} \!\! \pmod{p^3\ZZ_p} .
\]
Since $p(l-1)+1 = l + (p-1)(l-1) \equiv l \modp{p-1}$ and $(p,l) \in \IRR_1$,
only the first term $S_{p(l-1)+1}(p)$ does not vanish. Note that $p \notdiv p(l-1)+1$,
therefore we get
\[
   0 \equiv B_{l+(p-1)(l-1)} \equiv \BN(l+(p-1)(l-1)) \pmod{p^2}
\]
if and only if $(p,l,l-1) \in \IRP_2$.
\end{proofof}

\begin{remark}
From Propositions \ref{prop-delta-sn} and \ref{prop-delta-s1-s2},
the conditions $\Delta_{(p,\,l)} \neq 0$ and $(p,l,l-1) \notin \IRP_2$ are equivalent to the system
\begin{eqnarray*}
 &&  l \, S_{l+p-1}(p) - (l-1) \, S_{l}(p) \not\equiv 0  \pmod{p^3}  , \\
 &&  l \, S_{l+p-1}(p) - (l-2) \, S_{l}(p) \not\equiv 0  \pmod{p^3}  .
\end{eqnarray*}
\end{remark}

\section{The singular case}
\setcounter{equation}{0}

In Section \ref{sect-padic} we have derived most of the results
assuming the $\Delta$-Conjecture. Theorem \ref{theor-zeta-prod-psi}
conjecturally describes a closed formula for $\zeta(1-n)$ by zeros $\chi_{(p,\,l)}$.
The following theorem gives an equivalent formulation for the Bernoulli numbers.

\begin{theorem} \label{theor-bn-prod-psi}
Let $n$ be an even positive integer. Under the assumption of the $\Delta$-Conjecture
we have
\[
   B_n = (-1)^{\frac{n}{2}-1} \prod_{p-1 \notdiv n} |n|_p^{-1}
     \!\! \prod\limits_{(p,l)\in\IRR_1 \atop l \equiv n \,\mods{p-1}}
       \hspace*{-3ex} |p \, (\chi_{(p,\,l)} - {\textstyle \frac{n-l}{p-1}})|_p^{-1}
          \;\; \prod\limits_{p-1 \pdiv n} p^{-1} .
\]
\end{theorem}

\begin{proof}
We have to modify the formula of Theorem \ref{theor-zeta-prod-psi}.
The product formula gives
\begin{equation} \label{eqn-loc-zeta-prod-psi}
   1 = \prod_{p \, \in \, \PP \, \cup \, \{\infty\}} \hspace*{-2ex} |n|_p
     = |n|_\infty \prod_{p-1 \pdiv n} |n|_p \prod_{p-1 \notdiv n} |n|_p \,.
\end{equation}
Since $-B_n/n = \zeta(1-n)$, the proposed formula follows easily.
\end{proof}

To get an unconditional formula for $B_n$ resp.\ $\zeta(1-n)$
we have to include the case of a singular $\Delta_{(p,\,l)}$.
However, no such singular $\Delta_{(p,\,l)}$ has been found yet.
Theorem \ref{theor-delta-irr-0} describes the more complicated
behavior of related irregular pairs of higher order
in the singular case which can be described by a rooted $p$-ary tree,
see Diagram \ref{diag-tree}.

Let $(p,l) \in \IRR_1$ with a singular $\Delta_{(p,\,l)}$.
We construct the rooted $p$-ary tree of related irregular pairs of higher order
which is a consequence of Theorem \ref{theor-delta-irr-0}.
Each node contains one related irregular pair of higher order. Note that these pairs are not
necessarily distinct. We denote this tree by $T^0_{(p,\,l)}$ assigned to the root node $(p,l)$.

The tree $T^0_{(p,\,l)}$ has the property that each node of height $r$ lies in $\IRR_{r+1}$.
A tree $T^0_{(p,\,l)} = \{ (p,l) \}$ is called a \textsl{trivial tree} having height 0.
If the tree $T^0_{(p,\,l)}$ is of height $\geq 1$, then it
contains the root node $(p,l)$ and its $p$ child nodes
$(p,l+j\eulerphi(p)) \in \IRR_2$ for $j=0,\ldots,p-1$.

In the nonsingular case, we have a zero of the $p$-adic zeta function. In contrast
the singular case does not guarantee that related irregular pairs of higher order exist at all.
The discovery of a singular $\Delta_{(p,\,l)}$ is not incompatible with
Theorem \ref{theor-bn-prod-psi} but the formula becomes more complicated,
because we then have to consider the complete tree $T^0_{(p,\,l)}$.
By combining both cases we obtain an unconditional formula which is given
by the following theorem. Recall Definition \ref{def-ordn-irrpair}.

\begin{theorem} \label{theor-bn-prod-psi-0}
Let $n$ be an even positive integer, then
\begin{eqnarray*}
  B_n &=& (-1)^{\frac{n}{2}-1} \prod_{p-1 \notdiv n} |n|_p^{-1}
    \!\! \prod\limits_{(p,l)\in\IRR_1, \,\Delta_{(p,\,l)} \neq 0 \atop l \equiv n \,\mods{p-1}}
      \hspace*{-5ex} |p \, (\chi_{(p,\,l)} - {\textstyle \frac{n-l}{p-1}})|_p^{-1} \\
  && \quad \prod\limits_{(p,l)\in\IRR_1, \,\Delta_{(p,\,l)} = 0 \atop l \equiv n \,\mods{p-1}}
      \hspace*{-5ex} p^{1+h^0_{(p,\,l)}(n)}
        \;\; \prod\limits_{p-1 \pdiv n} p^{-1}
\end{eqnarray*}
with the height $h^0_{(p,\,l)}$ of $n$ defined by
\[
   h^0_{(p,\,l)}(n) = \max \left\{ \mathrm{height}( (p,l') ) :\,
     (p,l') \in T^0_{(p,\,l)} \cap \{ (p, n \modres \eulerphi(p^\nu) ) \}_{\nu \geq 1}
       \right\} .
\]
Moreover, $h^0_{(p,\,l)}(n) = 0$ $\iff$ the tree $T^0_{(p,\,l)}$ is trivial.
\end{theorem}

\begin{proof}
The case $\Delta_{(p,\,l)} \neq 0$ is already covered by Theorem \ref{theor-bn-prod-psi}.
Next we assume that $\Delta_{(p,\,l)} = 0$ with a given tree $T^0_{(p,\,l)}$
where $n \equiv l \modp{p-1}$.
As a consequence of the construction of $T^0_{(p,\,l)}$ and Remark \ref{rem-ordn-irrpair},
we have to determine the maximal height of a node $(p,l_{\nu,j}) \in T^0_{(p,\,l)} \cap \IRR_\nu$
where $(p,l_{\nu,j}) = ( p, n \modres \eulerphi(p^\nu))$.
The root node $(p,l)$ has height 0, so the exponent equals $1 + h^0_{(p,\,l)}(n)$.

If the tree $T^0_{(p,\,l)}$ has the height $\geq 1$, then
$(p, n \modres \eulerphi(p^2)) \in T^0_{(p,\,l)}$; this implies that
$h^0_{(p,\,l)}(n) \geq 1$. A trivial tree $T^0_{(p,\,l)}$ implies that
$h^0_{(p,\,l)}(n) = 0$. Conversely, if $h^0_{(p,\,l)}(n) = 0$, then
the height of the tree $T^0_{(p,\,l)}$ must be zero, otherwise we would
get a contradiction.
\end{proof}

\begin{corl}
Let $n$ be an even positive integer, then
\[
  \zeta(1-n) = (-1)^{\frac{n}{2}}
     \hspace*{-2ex} \prod\limits_{(p,l)\in\IRR_1, \,\Delta_{(p,\,l)} \neq 0
         \atop l \equiv n \,\mods{p-1}}
       \hspace*{-5ex} |p \, (\chi_{(p,\,l)} - {\textstyle \frac{n-l}{p-1}})|_p^{-1}
     \hspace*{-2ex} \prod\limits_{(p,l)\in\IRR_1, \,\Delta_{(p,\,l)} = 0
         \atop l \equiv n \,\mods{p-1}}
       \hspace*{-5ex} p^{1+h^0_{(p,\,l)}(n)}
     \; \prod\limits_{p-1 \pdiv n} \frac{|n|_p}{p}
\]
with $h^0_{(p,\,l)}$ as defined above.
\end{corl}

\begin{proof}
This is a reformulation of Theorem \ref{theor-bn-prod-psi-0} by
$\zeta(1-n) = -B_n/n$ and \refeqn{eqn-loc-zeta-prod-psi}.
\end{proof}

\section{An extension of Adams' theorem}
\setcounter{equation}{0}
\label{sect-adam}

\renewcommand{\thefootnote}{\fnsymbol{footnote}}

Let $n$ be an even positive integer.
The trivial factor of $B_n$, given by \refeqn{eqn-adams-triv-div}, is a consequence
of the implication, known as Adams' theorem, that
\[
   p^r \pdiv n \quad \mbox{with} \quad p-1 \notdiv n
     \qquad \Longrightarrow \qquad p^r \pdiv B_n
\]
where $p$ is a prime and $r$ is some positive integer.
It was, however, never proved by Adams. In 1878 he computed a table of
Bernoulli numbers $B_{2m}$ for $m \leq 62$.
On the basis of this table he conjectured%
\footnote[2]{
  ``... I have also observed that if $p$ be a prime factor of $n$ which is not
  likewise a factor of the denominator of the $n$th number of Bernoulli,
  then the numerator of that number will be divisible by $p$.
  I have not succeeded, however, in obtaining a general proof of this
  proposition, though I have no doubt of its truth.''}
that $p \pdiv n$ implies
$p \pdiv B_n$ for primes with $p-1 \notdiv n$, see \cite{adams78}.
Note that the property that $B_n/n$ is a $p$-integer for
$p-1 \notdiv n$ is needed to formulate the Kummer
congruences \refeqn{eqn-kummer-congr}. The case $r=1$ of these congruences
was proved by Kummer \cite{kummer51} earlier in 1851.
\renewcommand{\thefootnote}{\arabic{footnote}}
\smallskip

By Theorem \ref{theor-bn-prod-psi-0} and the definitions of $h^0_{(p,\,l)}$ and $\chi_{(p,\,l)}$,
we can state an extended version of Adams' theorem.
We introduce the following notation.
We write $p^r \pdiveq n$ when $p^r \pdiv n$ but $p^{r+1} \notdiv n$, i.e., $r=\ord_p n$.

\begin{theorem} \label{theor-ext-adam}
Let $n$ be an even positive integer. Let $p$ be a prime with $p^r \pdiveq n$,
$r \geq 1$, and $p-1 \notdiv n$. Let $l \equiv n \modp{p-1}$ with $0 < l < p-1$.
Then $p^{r+\delta} \pdiveq B_n$ with the following cases:
\begin{enumerate}
\item If $p$ is regular, then $\delta = 0$,
\item If $p$ is irregular with $(p,l)\notin\IRR_1$, then $\delta = 0$,
\item If $p$ is irregular with $(p,l)\in\IRR_1$, $\Delta_{(p,\,l)} \neq 0$,
      then $\delta = 1 + \ord_p \, (\chi_{(p,\,l)} - {\textstyle \frac{n-l}{p-1}})$,
\item If $p$ is irregular with $(p,l)\in\IRR_1$, $\Delta_{(p,\,l)} = 0$,
      then $\delta = 1 + h^0_{(p,\,l)}(n)$.
\end{enumerate}
Additionally, in case (3) resp.\ (4), if $(p,l,l) \notin \IRP_2$,
then $\delta = 1$, otherwise $\delta \geq 2$.
\end{theorem}
\pagebreak

\begin{proof}
We have to consider the formula of Theorem \ref{theor-bn-prod-psi-0}.
The first product yields $p^r \pdiv B_n$. Only the second resp.\ third product
can give additional $p$-factors. Therefore case (1) and (2) are given by definition.
We can now assume that $(p,l)\in\IRR_1$.
\smallskip

Case (3): A nonsingular $\Delta_{(p,\,l)}$ provides
\[
   \delta = \ord_p |p \, (\chi_{(p,\,l)} - {\textstyle \frac{n-l}{p-1}})|_p^{-1} =
     1 + \ord_p \, (\chi_{(p,\,l)} - {\textstyle \frac{n-l}{p-1}}) .
\]
By assumption $n=p^r n'$ with some integer $n'$. We have to evaluate
\[
   d = \ord_p \, (\chi_{(p,\,l)} - {\textstyle \frac{n-l}{p-1}}) =
     \ord_p \, (p \,\chi_{(p,\,l)} - \chi_{(p,\,l)} + l - p^r n') .
\]
Since $r \geq 1$, we $p$-adically obtain
\[
   (p,l,l) \in \IRP_2 \iff \chi_{(p,\,l)} = l + s_3\,p + \ldots \iff d \geq 1 .
\]
Conversely, $(p,l,l) \notin \IRP_2$ yields $d=0$.
\smallskip

Case (4): A singular $\Delta_{(p,\,l)}$ provides $\delta = 1 + h^0_{(p,\,l)}(n)$.
The definition of $T^0_{(p,\,l)}$ and Theorem \ref{theor-bn-prod-psi-0} show that
\[
   (p,l,l) \notin \IRP_2 \iff \mbox{the tree\ } T^0_{(p,\,l)} \mbox{\ is trivial}
     \iff h^0_{(p,\,l)}(n)=0 .
\]
Conversely, $(p,l,l) \in \IRP_2$ yields $h^0_{(p,\,l)}(n) \geq 1$.
\end{proof}

So far, no $(p,l,l) \in \IRP_2$ has been found.
The following corollary theoretically shows examples where $\delta$
is arbitrary large.

\begin{corl}
Assume that $(p,l,\ldots,l) \in \IRP_{r+1}$ exists with some $r\geq1$. Set $n=lp^r$.
Then we have $p^r \pdiveq n$ and $p^{2r+1} \pdiv B_n$, i.e., $\delta \geq r+1$.
\end{corl}

\begin{proof}
By Definition \ref{def-pair-padic}
$(p,n) \in \IRR_{r+1}$ is associated with $(p,l,\ldots,l) \in \IRP_{r+1}$,
since $n = lp^r = \sum_{\nu=1}^{r+1} l \eulerphi( p^{\nu-1} )$.
Thus, $p^{r+1} \pdiv B_n/n$ and finally $p^{2r+1} \pdiv B_n$.
\end{proof}

\begin{remark}
As mentioned above, Johnson \cite{johnson74} calculated
the now called irregular pairs $(p,s_1,s_2) \in \IRP_2$
of order two for $p < 8000$. He also proved that $(p,l,l) \notin \IRP_2$ resp.\
$s_1 \neq s_2$ in that range.
In a similar manner, the nonexistence of irregular pairs $(p,l,l-1)$ of order two
plays an important role in Iwasawa theory as seen in Section~\ref{sect-iwasawa}.
One may conjecture that no such special irregular pairs
$(p,l,l)$ and $(p,l,l-1)$ of order two exist. But there is still a long way
to prove such results, even to understand properly which role the zeros $\chi_{(p,\,l)}$ play.
Now, we have the relation
\[
   (p,l,l) \notin \IRP_2  \iff  p^2 \notdiv \BN(lp)
     \iff  p^3 \notdiv B_{lp} .
\]
Yamaguchi \cite{yamaguchi76} also verified by calculation that $p^3 \notdiv B_{lp}$
for all irregular pairs $(p,l)$ with $p < 5500$, noting
that this was conjectured earlier by Morishima in general.
The condition $p^3 \notdiv B_{lp}$ is related to the second
case of FLT, see \cite[Thm.~9.4, p.~174]{washington97}. See also
\cite[Cor.~8.23, p.~162]{washington97} for a different context.
Under the assumption of the conjecture of Kummer--Vandiver and that
no $(p,l,l) \in \IRP_2$ exists, the second case of FLT is true for the
exponent $p$. For details we refer to the references cited above.
\end{remark}

The converse of Adams' theorem does not hold, but one can state a somewhat different
result which deals with the common prime factors of numerators and denominators of
Bernoulli numbers with indices close to each other.

\begin{prop} \label{prop-bn-num-denom}
Let $\mathcal{S} = \{ 2,4,6,8,10,14 \}$ be the set of all even indices $m$ where the numerator
of $|B_m/m|$ equals 1. Write $B_n = \Bnum_n/\Bden_n$ with $(\Bnum_n,\Bden_n)=1$.
Let $k, n$ be even positive integers with $k \in \mathcal{S}$ and $n-k \geq 2$. Then
\[
    D=(\Bnum_n,\Bden_{n-k}) \quad \mbox{implies} \quad D \pdiv n .
\]
Moreover, if $D > 1$ then $D=p_1 \cdots p_r$ with some $r \geq 1$.
The primes $p_1, \ldots, p_r$ are pairwise different and $p_\nu \notdiv \Bden_k$,
$p_\nu \notdiv B_n/n$ for $\nu=1,\ldots,r$.
\end{prop}

\begin{proof}
Assume that $D > 1$. We then have $D=p_1 \cdots p_r$ with some $r \geq 1$,
since $\Bden_{n-k}$ is squarefree by \refeqn{eqn-clausen-von-staudt}.
Let $\nu \in \{ 1, \ldots, r \}$. Since $p_\nu \pdiv \Bnum_n$ and
$p_\nu \pdiv \Bden_{n-k}$, we have $p_\nu - 1 \notdiv n$ and
$p_\nu - 1 \pdiv n-k$.
From this we can deduce that $p_\nu-1 \notdiv k$ and consequently
that $p_\nu \notdiv \Bden_k$.
Next we assume that $p_\nu \notdiv n$ or $p_\nu \pdiv B_n/n$.
Note that $p_\nu \pdiv \Bnum_n$ and $p_\nu \notdiv n$ imply that $p_\nu \pdiv B_n/n$,
but not conversely.
We can use the Kummer congruences \refeqn{eqn-kummer-congr} to obtain that
\[
   0 \equiv \frac{B_n}{n} \equiv \frac{B_k}{k} \pmod{p_\nu} ,
\]
since $n \equiv k \modp{p_\nu-1}$. By the definition of the set
$\mathcal{S}$ we have
\begin{equation} \label{eqn-loc-bn-num-denom-1}
  \frac{B_k}{k} \not\equiv 0 \pmod{p_\nu}
\end{equation}
which yields a contradiction. This shows that $p_\nu \pdiv n$ and
$p_\nu \notdiv B_n/n$. Finally it follows that $D \pdiv n$.
\end{proof}

Now, the set $\mathcal{S}$ cannot be enlarged, because
\refeqn{eqn-loc-bn-num-denom-1}
does not hold in general for numerators having prime factors.
For example, let $p=691$ and $n=12+(p-1)=702$, then we have
$p \pdiv B_{12}/12$ and $D=(\Bnum_n,\Bden_{n-12})=pc \notdiv n$
with some $c \geq 1$. On the other hand, one trivially obtains for
$k \in \mathcal{S}$, $p$ prime with $p-1 \notdiv k$, and $n=kp$
infinitely many examples of $D > 1$. In the following proposition,
Proposition \ref{prop-bn-num-denom} plays a crucial role.
Recall the definition of $S_n(m)$.

\begin{prop} \label{prop-bn-sn}
Let $n, m$ be positive integers with even $n$. For $r=1,\,2$ we have
\[
   m^{r+1} \pdiv S_n(m) \quad \iff  \quad m^{r} \pdiv B_n .
\]
\end{prop}

\begin{proof}
We can assume that $m > 1$, since $m=1$ is trivial.
The case $n=2$ follows by $B_2 = \frac{1}{6}$
and that $m^2 \notdiv \frac{1}{6} m(m-1)(2m-1) = S_2(m)$ for $m > 1$.
For now we assume that $n \geq 4$. We have, see \refeqn{eqn-loc-delta-sn}, that
\begin{equation} \label{eqn-loc-bn-sn-1}
  S_n(m) = B_n \, m + \binom{n}{2} B_{n-2} \frac{m^3}{3} +
    \sum_{k=3}^n \binom{n}{k} B_{n-k} \frac{m^{k+1}}{k+1} .
\end{equation}
By \refeqn{eqn-clausen-von-staudt} and the cases
$B_0=1$ and $B_1 = -\frac12$
the denominator of all nonzero Bernoulli numbers is squarefree.
For each prime power factor $p^s \pdiveq m$
and $k$ where $B_{n-k} \neq 0$ $(2 \leq k \leq n)$ we have
\begin{equation} \label{eqn-loc-bn-sn-2}
   \ord_p \left( \binom{n}{k} B_{n-k} \frac{m^{k+1}}{k+1} \right)
     \geq s(k+1) - 1 - \ord_p (k+1) \geq \lambda \, s
\end{equation}
with the cases: (1) $\lambda=1$ for $k \geq 2, \, p \geq 2$,
(2) $\lambda=2$ for $k \geq 2, \, p \geq 5$, and (3) $\lambda=3$ for
$k \geq 4, \, p \geq 5$. The critical cases to consider are $p=2,3,5$
and $s=1$. Now, we are ready to evaluate \refeqn{eqn-loc-bn-sn-1}
$\mods{m^r}$ for $r=1,2$.
Write $B_n = \Bnum_n/\Bden_n$ with $(\Bnum_n,\Bden_n)=1$.
\smallskip

Case $r=1$: Assume that $(m,\Bden_n) > 1$.
By \refeqn{eqn-loc-bn-sn-2} (case $k \geq 2, \, p \geq 2$) we obtain
\[
   S_n(m) \equiv B_n \, m \equiv \frac{\Bnum_n}{\Bden_n} \,
     m \not\equiv 0 \pmod{m} .
\]
Therefore, $(m,\Bden_n) = 1$ must hold which implies $2 \notdiv m$, $3 \notdiv m$, and
$p \geq 5$. Hence, by \refeqn{eqn-loc-bn-sn-2}
(case $k \geq 2, \, p \geq 5$), we can write
$S_n(m) \equiv B_n \, m \modp{m^2}$. This yields
\begin{equation} \label{eqn-loc-bn-sn-3}
   m^2 \pdiv S_n(m) \iff m \pdiv B_n .
\end{equation}

Case $r=2$: We have $m \pdiv B_n$ and $(m,6)=1$, because
either $m^2 \pdiv B_n$ or $m^3 \pdiv S_n(m)$ is assumed.
The latter case implies $m^2 \pdiv S_n(m)$ and therefore
$m \pdiv B_n$ by \refeqn{eqn-loc-bn-sn-3}.
Since $|\Bnum_4|=1$, we can assume that $n \geq 6$.
We then have $B_{n-3} = 0$ and
we can apply \refeqn{eqn-loc-bn-sn-2}
(case $k \geq 4, \, p \geq 5$) to obtain
\begin{equation} \label{eqn-loc-bn-sn-4}
  S_n(m) \equiv B_n \, m + \frac{n(n-1)
    \Bnum_{n-2}}{6 \Bden_{n-2}} \, m^3 \pmod{m^3} .
\end{equation}
Our goal is to show that the second term of the right side
of \refeqn{eqn-loc-bn-sn-4} vanishes,
but the denominator $\Bden_{n-2}$ could possibly remove prime factors from $m$.
Proposition \ref{prop-bn-num-denom} asserts that
$(\Bnum_n, \Bden_{n-2}) \pdiv n$. We also have
$(m, \Bden_{n-2}) \pdiv n$ since $m \pdiv B_n$.
This means that the factor $n$ contains those primes which $\Bden_{n-2}$
possibly removes from $m$. Therefore the second term of
\refeqn{eqn-loc-bn-sn-4} vanishes $\mods{m^3}$. The rest follows
by $S_n(m) \equiv B_n \, m \equiv 0 \modp{m^3}$.
\end{proof}

One cannot improve the value $r$ in general. Choose $p=37$ and $l=37580$.
Since $(p,l) \in \IRR_3$ we have $p^3 \pdiv B_l$, but $p^4 \notdiv S_l(p)$
which was checked with \textbf{Mathematica}.
\smallskip

\begin{expenv}\

\begin{enumerate}
\item We have $B_{42} = 1520097643918070802691/1806$.
  Since the numerator $\Bnum_{42}$ is a large irregular prime,
  we obtain for $m>1$ that
\[
   m^2 \pdiv S_{42}(m) \quad \Longleftrightarrow \quad
     m = 1520097643918070802691 .
\]
\item We have
  $\Bnum_{50} = 5^2 \cdot 417202699 \cdot 47464429777438199$ and
  $\Bden_{48} = 2 \cdot 3 \cdot 5 \cdot 7 \cdot 13 \cdot 17$. Hence,
  for $m>1$ we have
\[
   m^3 \pdiv S_{50}(m) \quad \Longleftrightarrow \quad m = 5 .
\]
\end{enumerate}
\end{expenv}

\begin{appendix}

\section{Calculations}
\setcounter{equation}{0}

\begin{tablenv} \label{tab-bn} $B_n$ and $B_n/n$.

\begin{center}
\begin{scriptsize}
\begin{tabular}{c|*{12}{c}}
$n$ & 0 & 1 & 2 & 4 & 6 & 8 & 10 & 12 & 14 & 16 & 18 & 20 \\ \hline \\[-2pt]
 $B_n$ & 1 & $-\frac{1}{2}$ & $\frac{1}{6}$ &
$-\frac{1}{30}$ &
  $\frac{1}{42}$ & $-\frac{1}{30}$ & $\frac{5}{66}$ & $-\frac{691}{2730}$ &
  $\frac{7}{6}$ & $-\frac{3617}{510}$ & $\frac{43867}{798}$ &
  $-\frac{174611}{330}$ \\ \\ \hline \\[-2pt]
$\frac{B_n}{n}$ & & & $\frac{1}{12}$ &
$-\frac{1}{120}$ &
  $\frac{1}{252}$ & $-\frac{1}{240}$ & $\frac{1}{132}$ & $-\frac{691}{32760}$ &
  $\frac{1}{12}$ & $-\frac{3617}{8160}$ & $\frac{43867}{14364}$ &
  $-\frac{174611}{6600}$ \\ \\
\end{tabular}
\end{scriptsize}
\end{center}
\medskip

\end{tablenv}

\begin{tablenv} \label{tab-irrpair-ord-100}
Calculated irregular pairs of order 100 of primes 37, 59, and 67.
\smallskip

\begin{center}
\begin{scriptsize}
\begin{tabular}{*{11}{|r}|} \hline
\multicolumn{11}{|c|}{Case $p=37$.} \\
\multicolumn{11}{|c|}{Zeros of the sequence $(s_\nu)$
  occur at index 19 and 81.} \\ \hline\hline
 $s_\nu$ & 1 & 2 & 3 & 4 & 5 & 6 & 7 & 8 & 9 & 10 \\ \hline
 0  & 32 & 7 & 28 & 21 & 30 & 4 & 17 & 26 & 13 & 32 \\ \hline
 10 & 35 & 27 & 36 & 32 & 10 & 21 & 9 & 11 & 0 & 1 \\ \hline
 20 & 13 & 6 & 8 & 10 & 11 & 10 & 11 & 32 & 13 & 30 \\ \hline
 30 & 10 & 6 & 8 & 2 & 12 & 1 & 8 & 2 & 5 & 3 \\ \hline
 40 & 10 & 19 & 8 & 4 & 7 & 19 & 27 & 33 & 29 & 29 \\ \hline
 50 & 11 & 2 & 23 & 8 & 34 & 5 & 8 & 35 & 35 & 13 \\ \hline
 60 & 31 & 29 & 6 & 7 & 22 & 13 & 29 & 7 & 15 & 22 \\ \hline
 70 & 20 & 19 & 29 & 2 & 14 & 2 & 2 & 31 & 11 & 4 \\ \hline
 80 &  0 & 27 & 8 & 10 & 23 & 17 & 35 & 15 & 32 & 22 \\ \hline
 90 & 14 & 7 & 18 & 8 & 3 & 27 & 35 & 33 & 31 & 6 \\ \hline
\end{tabular}
\end{scriptsize}
\end{center}

\begin{center}
\begin{scriptsize}
\begin{tabular}{*{11}{|r}|} \hline
\multicolumn{11}{|c|}{Case $p=59$.} \\
\multicolumn{11}{|c|}{Zeros of the sequence $(s_\nu)$
  occur at index 31 and 95.} \\ \hline\hline
 $s_\nu$ & 1 & 2 & 3 & 4 & 5 & 6 & 7 & 8 & 9 & 10 \\ \hline
 0  & 44 & 15 & 25 & 40 & 36 & 18 & 11 & 17 & 28 & 58 \\ \hline
 10 &  9 & 51 & 13 & 25 & 41 & 44 & 17 & 43 & 35 & 21 \\ \hline
 20 & 10 & 21 & 38 &  9 & 12 & 40 & 43 & 45 & 30 & 41 \\ \hline
 30 &  0 &  3 & 25 & 34 & 49 & 45 &  9 & 19 & 48 & 57 \\ \hline
 40 & 11 & 13 & 29 & 28 & 44 & 41 & 37 & 33 & 29 & 43 \\ \hline
 50 &  8 & 57 & 12 & 48 & 15 & 15 & 53 & 57 & 16 & 51 \\ \hline
 60 & 16 & 54 & 30 &  9 & 26 &  8 & 49 & 22 & 58 & 11 \\ \hline
 70 & 42 & 28 & 36 & 33 & 45 & 24 & 32 & 18 & 12 & 29 \\ \hline
 80 & 45 & 40 & 27 & 19 & 40 & 41 & 11 & 42 & 49 & 35 \\ \hline
 90 & 41 & 57 & 54 & 33 &  0 & 34 & 34 & 49 &  6 & 31 \\ \hline
\end{tabular}
\end{scriptsize}
\end{center}

\begin{center}
\begin{scriptsize}
\begin{tabular}{*{11}{|r}|} \hline
\multicolumn{11}{|c|}{Case $p=67$.} \\
\multicolumn{11}{|c|}{Zeros of the sequence $(s_\nu)$
  occur at index 23 and 85.} \\ \hline\hline
 $s_\nu$ & 1 & 2 & 3 & 4 & 5 & 6 & 7 & 8 & 9 & 10 \\ \hline
 0  & 58 & 49 & 34 & 42 & 42 & 39 &  3 & 62 & 57 & 19 \\ \hline
 10 & 62 & 10 & 36 & 14 & 53 & 57 & 16 & 60 & 22 & 41 \\ \hline
 20 & 21 & 25 &  0 & 56 & 21 & 24 & 52 & 33 & 28 & 51 \\ \hline
 30 & 34 & 60 &  8 & 47 & 39 & 42 & 33 & 14 & 66 & 50 \\ \hline
 40 & 48 & 45 & 28 & 61 & 50 & 27 &  8 & 30 & 59 & 32 \\ \hline
 50 & 15 &  3 &  1 & 54 & 12 & 30 & 20 & 14 & 12 & 10 \\ \hline
 60 & 49 & 33 & 49 & 54 & 13 & 26 & 42 &  8 & 58 & 12 \\ \hline
 70 & 63 & 19 & 16 & 48 & 15 &  2 & 13 &  1 & 23 &  2 \\ \hline
 80 & 44 & 64 & 25 & 40 &  0 & 16 & 58 & 44 & 31 & 62 \\ \hline
 90 & 47 & 61 & 46 &  9 &  2 & 50 &  1 & 62 & 34 & 31 \\ \hline
\end{tabular}
\end{scriptsize}
\end{center}

\end{tablenv}
\pagebreak

\begin{tablenv} \label{tab-irrpair-ord-10}
Calculated irregular pairs of order 10 of primes below 1000.

\begin{center}
\begin{scriptsize}
\begin{tabular}{|c*{11}{|r}|} \hline
 $(p,l)$ & $\Delta_{(p,\,l)}$ & $s_1$ & $s_2$ & $s_3$ & $s_4$ & $s_5$ & $s_6$ &
    $s_7$ & $s_8$ & $s_9$ & $s_{10}$ \\ \hline \hline
(37,32) & 21 & 32 & 7 & 28 & 21 & 30 & 4 & 17 & 26 & 13 & 32 \\ \hline
(59,44) & 26 & 44 & 15 & 25 & 40 & 36 & 18 & 11 & 17 & 28 & 58 \\ \hline
(67,58) & 21 & 58 & 49 & 34 & 42 & 42 & 39 & 3 & 62 & 57 & 19 \\ \hline
(101,68) & 42 & 68 & 57 & 57 & 45 & 60 & 16 & 10 & 47 & 53 & 88 \\ \hline
(103,24) & 54 & 24 & 2 & 87 & 55 & 47 & 3 & 72 & 4 & 45 & 52 \\ \hline
(131,22) & 25 & 22 & 93 & 26 & 43 & 74 & 109 & 80 & 5 & 55 & 14 \\ \hline
(149,130) & 79 & 130 & 74 & 68 & 10 & 94 & 16 & 122 & 70 & 110 & 10 \\ \hline
(157,62) & 48 & 62 & 40 & 145 & 67 & 29 & 69 & 0 & 87 & 89 & 21 \\ \hline
(157,110) & 51 & 110 & 73 & 3 & 58 & 9 & 114 & 118 & 21 & 1 & 11 \\ \hline
(233,84) & 132 & 84 & 173 & 164 & 135 & 146 & 127 & 10 & 36 & 108 & 230 \\ \hline
(257,164) & 188 & 164 & 135 & 174 & 30 & 203 & 161 & 193 & 142 & 68 & 126 \\ \hline
(263,100) & 87 & 100 & 198 & 139 & 151 & 106 & 202 & 99 & 202 & 251 & 163 \\ \hline
(271,84) & 179 & 84 & 5 & 14 & 239 & 8 & 233 & 43 & 28 & 57 & 170 \\ \hline
(283,20) & 15 & 20 & 265 & 115 & 171 & 137 & 251 & 118 & 132 & 246 & 265 \\ \hline
(293,156) & 93 & 156 & 230 & 75 & 289 & 47 & 247 & 98 & 100 & 141 & 27 \\ \hline
(307,88) & 205 & 88 & 70 & 234 & 51 & 173 & 104 & 140 & 140 & 107 & 201 \\ \hline
(311,292) & 277 & 292 & 204 & 183 & 9 & 260 & 183 & 214 & 254 & 2 & 151 \\ \hline
(347,280) & 106 & 280 & 113 & 250 & 150 & 307 & 264 & 145 & 177 & 101 & 156 \\ \hline
(353,186) & 301 & 186 & 190 & 147 & 13 & 34 & 171 & 106 & 304 & 190 & 102 \\ \hline
(353,300) & 161 & 300 & 181 & 300 & 314 & 327 & 67 & 26 & 113 & 18 & 336 \\ \hline
(379,100) & 276 & 100 & 242 & 277 & 88 & 236 & 225 & 22 & 221 & 54 & 26 \\ \hline
(379,174) & 82 & 174 & 364 & 216 & 20 & 128 & 277 & 134 & 257 & 164 & 31 \\ \hline
(389,200) & 48 & 200 & 354 & 33 & 371 & 189 & 29 & 219 & 44 & 11 & 319 \\ \hline
(401,382) & 376 & 382 & 263 & 126 & 213 & 197 & 170 & 320 & 107 & 297 & 331 \\ \hline
(409,126) & 180 & 126 & 389 & 343 & 247 & 322 & 24 & 187 & 75 & 91 & 179 \\ \hline
(421,240) & 396 & 240 & 351 & 141 & 36 & 169 & 124 & 164 & 342 & 365 & 156 \\ \hline
(433,366) & 284 & 366 & 406 & 342 & 372 & 234 & 21 & 328 & 346 & 279 & 155 \\ \hline
(461,196) & 281 & 196 & 423 & 121 & 233 & 61 & 353 & 421 & 414 & 350 & 92 \\ \hline
(463,130) & 78 & 130 & 376 & 404 & 124 & 420 & 63 & 438 & 185 & 124 & 18 \\ \hline
(467,94) & 118 & 94 & 219 & 393 & 264 & 70 & 75 & 254 & 361 & 332 & 157 \\ \hline
(467,194) & 269 & 194 & 283 & 329 & 154 & 419 & 170 & 152 & 78 & 304 & 326 \\ \hline
(491,292) & 456 & 292 & 218 & 299 & 225 & 362 & 461 & 37 & 65 & 203 & 228 \\ \hline
(491,336) & 103 & 336 & 260 & 15 & 41 & 381 & 66 & 376 & 391 & 209 & 305 \\ \hline
(491,338) & 475 & 338 & 59 & 160 & 106 & 105 & 33 & 346 & 158 & 314 & 233 \\ \hline
(523,400) & 497 & 400 & 36 & 230 & 180 & 431 & 235 & 114 & 104 & 152 & 399 \\ \hline
(541,86) & 211 & 86 & 436 & 29 & 482 & 424 & 74 & 212 & 259 & 419 & 287 \\ \hline
(547,270) & 348 & 270 & 458 & 536 & 35 & 521 & 413 & 88 & 545 & 44 & 537 \\ \hline
(547,486) & 139 & 486 & 100 & 4 & 33 & 153 & 282 & 467 & 233 & 482 & 17 \\ \hline
(557,222) & 153 & 222 & 549 & 505 & 399 & 472 & 49 & 20 & 81 & 279 & 513 \\ \hline
(577,52) & 452 & 52 & 309 & 416 & 274 & 56 & 20 & 476 & 164 & 309 & 19 \\ \hline
(587,90) & 286 & 90 & 109 & 344 & 244 & 53 & 93 & 454 & 292 & 291 & 547 \\ \hline
(587,92) & 319 & 92 & 213 & 332 & 470 & 36 & 479 & 508 & 134 & 323 & 275 \\ \hline
(593,22) & 331 & 22 & 188 & 388 & 541 & 576 & 371 & 26 & 586 & 40 & 514 \\ \hline
(607,592) & 435 & 592 & 369 & 428 & 162 & 503 & 358 & 484 & 411 & 67 & 267 \\ \hline
(613,522) & 57 & 522 & 549 & 451 & 318 & 312 & 243 & 38 & 265 & 552 & 215 \\ \hline
(617,20) & 289 & 20 & 384 & 107 & 161 & 281 & 358 & 64 & 604 & 336 & 326 \\ \hline
(617,174) & 317 & 174 & 546 & 83 & 114 & 484 & 121 & 229 & 335 & 597 & 570 \\ \hline
(617,338) & 312 & 338 & 419 & 570 & 496 & 63 & 247 & 46 & 604 & 464 & 134 \\ \hline
(619,428) & 121 & 428 & 457 & 363 & 526 & 36 & 179 & 79 & 170 & 485 & 47 \\ \hline
(631,80) & 139 & 80 & 146 & 468 & 175 & 34 & 249 & 169 & 26 & 498 & 528 \\ \hline
(631,226) & 221 & 226 & 338 & 510 & 318 & 581 & 572 & 363 & 422 & 111 & 405 \\ \hline
(647,236) & 318 & 236 & 480 & 525 & 205 & 103 & 205 & 620 & 394 & 553 & 25 \\ \hline
(647,242) & 94 & 242 & 487 & 519 & 49 & 109 & 373 & 451 & 586 & 250 & 57 \\ \hline
(647,554) & 209 & 554 & 558 & 568 & 174 & 579 & 545 & 5 & 377 & 242 & 81 \\ \hline
(653,48) & 363 & 48 & 154 & 558 & 439 & 300 & 59 & 541 & 242 & 205 & 47 \\ \hline
(659,224) & 200 & 224 & 140 & 131 & 396 & 158 & 367 & 79 & 256 & 620 & 615 \\ \hline
(673,408) & 325 & 408 & 26 & 64 & 257 & 158 & 213 & 430 & 659 & 144 & 600 \\ \hline
(673,502) & 585 & 502 & 293 & 198 & 436 & 506 & 441 & 27 & 89 & 416 & 407 \\ \hline
\end{tabular}
\end{scriptsize}
\end{center}
\end{tablenv}

\begin{center}
\begin{scriptsize}
\begin{tabular}{|c*{11}{|r}|} \hline
 $(p,l)$ & $\Delta_{(p,\,l)}$ & $s_1$ & $s_2$ & $s_3$ & $s_4$ & $s_5$ & $s_6$ &
    $s_7$ & $s_8$ & $s_9$ & $s_{10}$ \\ \hline \hline
(677,628) & 440 & 628 & 504 & 457 & 324 & 461 & 88 & 532 & 653 & 89 & 244 \\ \hline
(683,32) & 477 & 32 & 266 & 20 & 625 & 119 & 190 & 13 & 190 & 222 & 214 \\ \hline
(691,12) & 611 & 12 & 496 & 104 & 197 & 607 & 590 & 303 & 96 & 461 & 152 \\ \hline
(691,200) & 592 & 200 & 496 & 333 & 578 & 93 & 160 & 436 & 611 & 215 & 278 \\ \hline
(727,378) & 398 & 378 & 683 & 722 & 169 & 391 & 150 & 694 & 210 & 228 & 130 \\ \hline
(751,290) & 164 & 290 & 481 & 37 & 181 & 27 & 31 & 71 & 8 & 36 & 164 \\ \hline
(757,514) & 554 & 514 & 364 & 164 & 375 & 7 & 720 & 750 & 273 & 592 & 643 \\ \hline
(761,260) & 462 & 260 & 729 & 680 & 274 & 188 & 464 & 183 & 283 & 52 & 235 \\ \hline
(773,732) & 517 & 732 & 147 & 306 & 278 & 370 & 412 & 89 & 340 & 637 & 223 \\ \hline
(797,220) & 375 & 220 & 369 & 279 & 501 & 300 & 168 & 530 & 534 & 747 & 268 \\ \hline
(809,330) & 88 & 330 & 52 & 743 & 100 & 336 & 157 & 759 & 348 & 43 & 736 \\ \hline
(809,628) & 18 & 628 & 773 & 629 & 623 & 160 & 494 & 339 & 244 & 463 & 274 \\ \hline
(811,544) & 381 & 544 & 424 & 100 & 346 & 749 & 624 & 220 & 410 & 313 & 62 \\ \hline
(821,744) & 704 & 744 & 621 & 319 & 498 & 427 & 50 & 21 & 237 & 305 & 809 \\ \hline
(827,102) & 105 & 102 & 164 & 443 & 469 & 568 & 671 & 183 & 372 & 512 & 464 \\ \hline
(839,66) & 269 & 66 & 135 & 305 & 36 & 40 & 659 & 431 & 326 & 591 & 293 \\ \hline
(877,868) & 480 & 868 & 554 & 279 & 714 & 821 & 520 & 76 & 565 & 104 & 22 \\ \hline
(881,162) & 789 & 162 & 372 & 330 & 89 & 244 & 27 & 229 & 418 & 438 & 89 \\ \hline
(887,418) & 611 & 418 & 76 & 698 & 835 & 872 & 130 & 319 & 217 & 439 & 573 \\ \hline
(929,520) & 607 & 520 & 433 & 27 & 711 & 366 & 902 & 838 & 7 & 351 & 805 \\ \hline
(929,820) & 706 & 820 & 749 & 156 & 59 & 913 & 480 & 432 & 114 & 129 & 491 \\ \hline
(953,156) & 24 & 156 & 720 & 516 & 620 & 229 & 251 & 77 & 805 & 689 & 477 \\ \hline
(971,166) & 715 & 166 & 538 & 594 & 897 & 509 & 355 & 749 & 180 & 174 & 96 \\ \hline
\end{tabular}
\end{scriptsize}
\end{center}
\smallskip

\begin{tablenv} \label{tab-calc-37-1}
Calculation: $n=1$, $p=37$, $l = 32$, $(p,l) \in \IRR_1$.

\begin{center}
\begin{scriptsize}
\begin{tabular}{|c|r|c|r|r|r|} \hline
$j$ & Index & $\alpha_j \modp{p^3}$ & $\equiv \modp{p^3}$
    & $\Delta_{\alpha_j} \modp{p^3}$ & $\Delta_{\alpha_j} \modp{p^2}$ \\ \hline\hline
0 &  32 & 3941/2720 & 42144 & 45827 & 650 \\ \hline
1 &  68 & 2587/15   & 37318 & 49934 & 650 \\ \hline
2 & 104 & 3821/1272 & 36599 & 30768 & 650 \\ \hline
3 & 140 & 6497/7198 & 16714 & \multicolumn{2}{r|}{$\Delta_{(p,\,l)} = 21$} \\ \hline
\end{tabular}
\end{scriptsize}
\end{center}

Using Proposition \ref{prop-irrpair-rn} with $r=3$ and $(r-1)n=2$
yields $s \equiv 1043 \modp{p^2}$ and $l_3 = 32 + s \eulerphi(p) = 37580$.
We obtain $(37,284) \in \IRR_2$, $(37,37580) \in \IRR_3$, and
$(37,32,7,28) \in \IRP_3$.
\end{tablenv}
\smallskip

\begin{tablenv} \label{tab-calc-37-2}
Calculation: $n=3$, $p=37$, $l = 37580$, $(p,l) \in \IRR_3$.

\begin{center}
\begin{scriptsize}
\begin{tabular}{|c|r|c|r|r|} \hline
$j$ & Index & $\alpha_j \modp{p^3}$ & $\equiv \modp{p^3}$
    & $\Delta_{\alpha_j} \modp{p^3}$ \\ \hline\hline
0 &  37580 & 11241/22913 & 24645 & 45827  \\ \hline
1 &  86864 & 49609/46188 & 19819 & 45827  \\ \hline
2 & 136148 & 5261/24     & 14993 & $\Delta_{(p,\,l)} = 21$ \\ \hline
\end{tabular}
\end{scriptsize}
\end{center}

Using Proposition \ref{prop-irrpair-2n} yields $s \equiv 6607 \modp{p^3}$
and $l_6 = 37580 + s \eulerphi(p^3) = 325656968$. We obtain
$(37,325656968) \in \IRR_6$, $(37,55777784) \in \IRR_5$,
$(37,1072544) \in \IRR_4$, and $(37,32,7,28,21,30,4) \in \IRP_6$.
\end{tablenv}
\smallskip

\begin{tablenv} \label{tab-calc-37-3}
Calculation: $n=3$, $p=37$, $l = 37580$, $(p,l) \in \IRR_3$.

\begin{center}
\begin{scriptsize}
\begin{tabular}{|c|r|c|r|} \hline
$j$ & Index & $\alpha_j \modp{p^9}$ & $\equiv \modp{p^9}$ \\ \hline\hline
0 &  37580 & 3791602112159/3307480      & 45520991695194 \\ \hline
1 &  86864 & 1046892158059/484258896735 & 47985230204445 \\ \hline
2 & 136148 & 13280633201029/15          & 70198303437443 \\ \hline
3 & 185432 & 8822143378793/98280020     & 73479320052104 \\ \hline
\end{tabular}
\end{scriptsize}
\end{center}

Using Proposition \ref{prop-calc-irrpair} with $r=4$ and $(r-1)n=9$ yields
the sequence $21, 30, 4, \ldots, 27$ which provides
$(37,32,7,28,21,30,4,17,26,13,32,35,27) \in \IRP_{12}$.
\end{tablenv}
\pagebreak

Note that Tables \ref{tab-irrpair-ord-100} and \ref{tab-irrpair-ord-10}
were calculated with smallest possible indices of the Bernoulli numbers using
Proposition \ref{prop-calc-irrpair}; they agree with these results above.
Additionally, the results were checked by Corollary \ref{corl-zeta-chi-z}
and Proposition \ref{prop-calc-chi}.
The program \mbox{\textbf{calcbn}} \cite[Section 2.7]{kellner02}
was used to calculate these large Bernoulli numbers extremely quickly.

\end{appendix}

\addcontentsline{toc}{section}{References}

\bibliographystyle{plain}

\medskip

\textsc{Reitstallstr.~7, 37073 G\"ottingen, Germany}

\textsl{E-mail address}: \textbf{bk@bernoulli.org}

\end{document}